%% file: LiuZhao09TSP_sub2_revision.tex
\documentclass[11pt,onecolumn,journal,draftclsnofoot,doublespace]{IEEEtran}
\usepackage{times,amsmath,amssymb,graphicx,epsf,psfrag,epsfig,cite}

\input macros

\def\ie{{\it i.e.,\ \/}}

\def\defeq{{\,\stackrel{\Delta}{=}}\,}

\def\scalefig#1{\epsfxsize #1\textwidth}
\def\nn{{\nonumber}}

\newtheorem{lemma}{Lemma}
\newtheorem{theorem}{Theorem}

\newtheorem{definition}{Definition}

\begin{document}
\title{Distributed Learning in Multi-Armed Bandit with Multiple Players}
\author{ Keqin Liu, ~~~ Qing Zhao\\
 University of California, Davis, CA 95616\\
 \{kqliu, qzhao\}@ucdavis.edu}
\date{}
\markboth{Submitted to {\em IEEE Transactions on Signal Processing},
December, 2009.}{Liu and Zhao} \maketitle

\vspace{-2em}

\begin{abstract}
We \footnotetext{This work was supported by the Army Research
Laboratory NS-CTA under Grant W911NF-09-2-0053, by the Army Research
Office under Grant W911NF-08-1-0467, and by the National Science
Foundation under Grants CNS-0627090 and CCF-0830685. Part of this
result will be presented at ICASSP 2010.} formulate and study a
decentralized multi-armed bandit (MAB) problem. There are $M$
distributed players competing for $N$ independent arms. Each arm,
when played, offers i.i.d. reward according to a distribution with
an unknown parameter. At each time, each player chooses one arm to
play without exchanging observations or any information with other
players. Players choosing the same arm collide, and, depending on
the collision model, either no one receives reward or the colliding
players share the reward in an arbitrary way. We show that the
minimum system regret of the decentralized MAB grows with time at
the same logarithmic order as in the centralized counterpart where
players act collectively as a single entity by exchanging
observations and making decisions jointly. A decentralized policy is
constructed to achieve this optimal order while ensuring fairness
among players and without assuming any pre-agreement or information
exchange among players. Based on a Time Division Fair Sharing (TDFS)
of the $M$ best arms, the proposed policy is constructed and its
order optimality is proven under a general reward model.
Furthermore, the basic structure of the TDFS policy can be used with
any order-optimal single-player policy to achieve order optimality
in the decentralized setting. We also establish a lower bound on the
system regret growth rate for a general class of decentralized
polices, to which the proposed policy belongs. This problem finds
potential applications in cognitive radio networks, multi-channel
communication systems, multi-agent systems, web search and
advertising, and social networks.
\end{abstract}

\vspace{0.5em}

\begin{IEEEkeywords}
Decentralized multi-armed bandit, system regret, distributed
learning, cognitive radio, web search and advertising, multi-agent
systems.
\end{IEEEkeywords}

\section{Introduction}\label{sec:intro}

\subsection{The Classic MAB with A Single Player}
In the classic multi-armed bandit (MAB) problem, there are $N$
independent arms and a single player. Playing arm $i$ ($i=1,\cdots,
N$) yields i.i.d. random rewards with a distribution parameterized
by an unknown $\theta_i$. At each time, the player chooses one arm
to play, aiming to maximize the total expected reward in the long
run. Had the reward model of each arm been known to the player, the
player would have always chosen the arm that offers the maximum
expected reward. Under an unknown reward model, the essence of the
problem lies in the well-known tradeoff between exploitation (aiming
at gaining immediate reward by choosing the arm that is suggested to
be the best by past observations) and exploration (aiming at
learning the unknown reward model of all arms to minimize the
mistake of choosing an inferior arm in the future).

Under a non-Bayesian formulation where the unknown parameters
$\Theta\defeq(\theta_1,\cdots,\theta_N)$ are considered
deterministic, a commonly used performance measure of an arm
selection policy is the so-called \emph{regret} or \emph{the cost of
learning} defined as the reward loss with respect to the case with
known reward models. Since the best arm can never be perfectly
identified from finite observations (except certain trivial
scenarios), the player can never stop learning and will always make
mistakes. Consequently, the regret of any policy grows with time.

An interesting question posed by Lai and Robbins in
1985~\cite{Lai&Robbins85AAM} is on the minimum rate at which the
regret grows with time. They showed in~\cite{Lai&Robbins85AAM} that
the minimum regret grows at a logarithmic order under certain
regularity conditions. The best leading constant was also obtained,
and an optimal policy was constructed under a general reward model
to achieve the minimum regret growth rate (both the logarithmic
order and the best leading constant). In 1987, Anantharam \etal
extended Lai and Robbins's results to MAB with multiple plays:
exactly $M~(M<N)$ arms can be played simultaneously at each
time~\cite{Anantharam&etal87TAC}. They showed that allowing multiple
plays changes only the leading constant but not the logarithmic
order of the regret growth rate. They also extended Lai-Robbins
policy to achieve the optimal regret growth rate under multiple
plays.

\subsection{Decentralized MAB with Distributed Multiple Players}
In this paper, we formulate and study a decentralized version of the
classic MAB, where we consider $M$ ($M<N$) distributed players. At
each time, a player chooses one arm to play based on its {\em local}
observation and decision history. Players do not exchange
information on their decisions and observations. Collisions occur
when multiple players choose the same arm, and, depending on the
collision model, either no one receives reward or the colliding
players share the reward in an arbitrary way. The objective is to
design distributed policies for each player in order to minimize,
under any unknown parameter set $\Theta$, the rate at which the
system regret grows with time. Here the system regret is defined as
the reward loss with respect to the maximum system reward obtained
under a known reward model and with centralized scheduling of the
players.

The single-player MAB with multiple plays considered
in~\cite{Anantharam&etal87TAC} is equivalent to a \emph{centralized}
MAB with multiple players. If all $M$ players can exchange
observations and make decisions jointly, they act collectively as a
single player who has the freedom of choosing $M$ arms
simultaneously. As a direct consequence, the logarithmic order of
the minimum regret growth rate established
in~\cite{Anantharam&etal87TAC} provides a lower bound on the minimum
regret growth rate in a \emph{decentralized} MAB where players
cannot exchange observations and must make decisions independently
based on their individual local observations\footnote{While
intuitive, this equivalence requires the condition that the $M$ best
arms have nonnegative means (see Sec.~\ref{sec:orderoptimal} for
details).}.

\subsection{Main Results}
In this paper, we show that in a decentralized MAB where players can
only learn from their individual observations and collisions are
bound to happen, the system can achieve the same logarithmic order
of the regret growth rate as in the centralized case. Furthermore,
we show that this optimal order can be achieved under a fairness
constraint that requires all players accrue reward at the same rate.

A decentralized policy is constructed to achieve the optimal order
under the fairness constraint. The proposed policy is based on a
Time Division Fair Sharing (TDFS) of the $M$ best arms where no
pre-agreement on the time sharing schedule is needed. It is
constructed and its order optimality is proven under a general
reward model. The TDFS decentralized policy thus applies to and
maintains its order optimality in a wide range of potential
applications as detailed in Sec.~\ref{subsec:app}. Furthermore, the
basic structure of TDFS is not tied with a specific single-player
policy. It can be used with any single-player policy to achieve an
efficient and fair sharing of the $M$ best arms among $M$
distributed players. More specifically, if the single-player policy
achieves the optimal logarithmic order in the centralized setting,
then the corresponding TDFS policy achieves the optimal logarithmic
order in the decentralized setting. The order optimality of the TDFS
policy is also preserved when player's local polices are built upon
different order-optimal single-player polices.

We also establish a lower bound on the leading constant in the
regret growth rate for a general class of decentralized polices, to
which the proposed TDFS policy belongs. This lower bound is tighter
than the trivial bound provided by the centralized MAB considered
in~\cite{Anantharam&etal87TAC}, which indicates, as one would have
expected, that decentralized MAB is likely to incur a larger leading
constant than its centralized counterpart.

\subsection{Applications}\label{subsec:app}
With its general reward model, the TDFS policy finds a wide area of
potential applications. We give a few examples below.

Consider first a cognitive radio network where secondary users
independently search for idle channels temporarily unused by primary
users~\cite{Zhao&Sadler:07SPM}. Assume that the primary system
adopts a slotted transmission structure and the state (busy/idle) of
each channel can be modeled by an i.i.d. Bernoulli process. At the
beginning of each slot, multiple distributed secondary users need to
decide which channel to sense (and subsequently transmit if the
chosen channel is idle) without knowing the channel occupancy
statistics (\ie the mean $\theta_i$ of the Bernoulli process). An
idle channel offers one unit of reward. When multiple secondary
users choose the same channel, none or only one receives reward
depending on whether carrier sensing is implemented.

Another potential application is opportunistic transmission over
wireless fading channels~\cite{Liu&etal:01JSAC}. In each slot, each
user senses the fading realization of a selected channel and chooses
its transmission power or date rate accordingly. The reward can
model energy efficiency (for fixed-rate transmission) or throughput.
The objective is to design distributed channel selection policies
under unknown fading statistics.

Consider next a multi-agent system, where each agent is assigned to
collect targets among multiple locations (for example, ore mining).
When multiple agents choose the same location, they share the reward
according to a certain rule. The log-Gaussian distribution with an
unknown mean may be used as the reward model when the target is
fish~\cite{Myers&Pepin:90Biometrics} or ore~\cite{Sharp:76EG}.

Another potential application is Internet advertising where multiple
competing products select Web sites to posts advertisements. The
reward obtained from the selected Web site is measured by the number
of interested viewers whose distribution can be modeled as
Poisson~\cite{Bean&09book}.

\subsection{Related Work}
In the context of the classic MAB, there have been several attempts
at developing index-type policies that are simpler than Lai-Robbins
policy by using a single sample mean based statistic [8,9]. However,
such a simpler form of the policy was obtained at the price of a
more restrictive reward model and/or a larger leading constant in
the logarithmic order. In particular, the index policy proposed
in~\cite{Agrawal95AAP} mainly targets at several members of the
exponential family of reward distributions. Auer \etal proposed
in~\cite{Auer&etal02ML} several index-type policies that achieve
order optimality for reward distributions with a known finite
support.

Under a Bernoulli reward model in the context of cognitive radio,
MAB with multiple players was considered
in~\cite{Lai&etal08Asilomar} and \cite{Anandkumar&10INFOCOM}.
In~\cite{Lai&etal08Asilomar}, a heuristic policy based on histogram
estimation of the unknown parameters was proposed. This policy
provides a linear order of the system regret rate, thus cannot
achieve the maximum average reward. In~\cite{Anandkumar&10INFOCOM},
Anandkumar \etal have independently established order-optimal
distributed policies by extending the index-type single-user
policies proposed in~\cite{Auer&etal02ML}. Compared
to~\cite{Anandkumar&10INFOCOM}, the TDFS policy proposed here
applies to more general reward models (for example, Gaussian and
Poisson reward distributions that have infinite support). It thus
has a wider range of potential applications as discussed in
Sec.~\ref{subsec:app}. Furthermore, the policies proposed
in~\cite{Anandkumar&10INFOCOM} are specific to the single-player
polices proposed in~\cite{Auer&etal02ML}, whereas the TDFS policy
can be used with any order-optimal single-player policy to achieve
order optimality in the decentralized setting. Another difference
between the policies proposed in~\cite{Anandkumar&10INFOCOM} and the
TDFS policy is on user fairness. The policies
in~\cite{Anandkumar&10INFOCOM} orthogonalize users into different
channels that offer different throughput, whereas the TDFS policy
ensures that each player achieves the same time-average reward at
the same rate. One policy given in~\cite{Anandkumar&10INFOCOM} does
offer probabilistic fairness in the sense that all users have the
same chance of settling in the best channel. However, given that the
policy operates over an infinite horizon and the order optimality is
asymptotic, each user only sees one realization in its lifetime,
leading to different throughput among users. A lower bound on the
achievable growth rate of the system regret is also given
in~\cite{Anandkumar&10INFOCOM}, which is identical to the bound
developed here. The derivation of the lower bound in this work,
however, applies to a more general class of policies and, first
given in~\cite{Liu&Zhao09ArXiv}, precedes that
in~\cite{Anandkumar&10INFOCOM}. In terms of using collision history
to orthogonalizing players without pre-agreement, the basic idea
used in this work is similar to that in~\cite{Anandkumar&10INFOCOM}.
The difference is that in the TDFS policy, the players are
orthogonalized via settling at different offsets in their time
sharing schedule, while in~\cite{Anandkumar&10INFOCOM}, players are
orthogonalized to different channels. Recently, results obtained
in~\cite{Anandkumar&10INFOCOM} have been extended to incorporate
unknown number of users for Bernoulli reward model in the context of
cognitive radio~\cite{Anandkumar&09JSAC_sub}.

A variation of centralized MAB in the context of cognitive radio has
been considered in~\cite{Gai&etal:10DySPANsub} where a channel
offers independent Bernoulli rewards with different (unknown) means
for different users. This more general model captures contention
among secondary users that are in the communication range of
different sets of primary users. A centralized policy that assumes
full information exchange and cooperation among users is proposed
which achieves the logarithmic order of the regret growth rate. We
point out that the logarithmic order of the TDFS policy is preserved
in the decentralized setting when we allow players to experience
different reward distributions on the same arm provided that players
have the common set of the $M$ best arms and each of the $M$ best
arms has the same mean across players. It would be interesting to
investigate whether the same holds when players experience different
means on each arm.

\subsection{Notations}
Let $|\Ac|$ denote the cardinality of set $\Ac$. For two sets $\Ac$
and $\Bc$, let $\Ac\backslash\Bc$ denote the set consisting of all
elements in $\Ac$ that do not belong to $\Bc$. For two positive
integers $k$ and $l$, define $k\oslash l\defeq
((k-1)~\mbox{mod}~l)+1$, which is an integer taking values from
$1,2,\cdots,l$. Let $\Pr_{\theta}\{\cdot\}$ denote the probability
measure when the unknown parameter in the associated distribution
equals to $\theta$.

\section{Classic Results on Single-Player MAB}\label{sec:classic}
In this section, we give a brief review of the main results
established
in~\cite{Lai&Robbins85AAM,Anantharam&etal87TAC,Agrawal95AAP,Auer&etal02ML}
on the classic MAB with a single player.

\subsection{System Regret}
Consider an $N$-arm bandit with a single player. At each time $t$,
the player chooses exactly $M$ ($1\le M<N$) arms to play. Playing
arm $i$ yields i.i.d. random reward $S_i(t)$ drawn from a univariate
density function $f(s;\theta_i)$ parameterized by $\theta_i$. The
parameter set $\Theta\defeq(\theta_1,\cdots,\theta_N)$ is unknown to
the player. Let $\mu(\theta_i)$ denote the mean of $S_i(t)$ under
the density function $f(s;\theta_i)$. Let
$I(\theta,\theta')=\mathbb{E}_{\theta}[\log[\frac{f(y,\theta)}{f(y,\theta')}]
]$ be the Kullback-Liebler distance that measures the dissimilarity
between two distributions parameterized by $\theta$ and $\theta'$,
respectively.

An arm selection policy $\pi=\{\pi(t)\}_{t=1}^{\infty}$ is a series
of functions, where $\pi(t)$ maps the previous observations of
rewards to the current action that specifies the set of $M$ arms to
play at time $t$. The system performance under policy $\pi$ is
measured by the system regret $R_T^{\pi}(\Theta)$ defined as the
expected total reward loss up to time $T$ under policy $\pi$
compared to the ideal scenario that $\Theta$ is known to the player
(thus the $M$ best arms are played at each time). Let $\sigma$ be a
permutation of $\{1,\cdots,N\}$ such that
$\mu(\theta_{\sigma(1)})\ge\mu(\theta_{\sigma(2)})\ge\cdots\ge\mu(\theta_{\sigma(N)})$.
We have
\[R_T^{\pi}(\Theta)\defeq T\Sigma_{j=1}^M\mu(\theta_{\sigma(j)})-\mathbb{E}_{\pi}[\Sigma_{t=1}^TS_{\pi(t)}(t)],\]
where $S_{\pi(t)}(t)$ is the random reward obtained at time $t$
under action $\pi(t)$, and $\mathbb{E}_{\pi}[\cdot]$ denotes the
expectation with respect to policy $\pi$. The objective is to
minimize the rate at which $R_T(\Theta)$ grows with $T$ under any
parameter set $\Theta$ by choosing an optimal policy $\pi^*$.

We point out that the system regret rate is a finer performance
measure than the long-term average reward. All policies with a
sublinear regret rate would achieve the maximum long-term average
reward. However, the difference in their performance measured by the
total expected reward accrued over a time horizon of length $T$ can
be arbitrarily large as $T$ increases. It is thus of great interest
to characterize the minimum regret rate and construct policies
optimal under this finer performance measure.

A policy is called {\em uniformly good} if for any parameter set
$\Theta$, we have $R_T^{\pi}(\Theta)=o(T^b)$ for any $b>0$. Note
that a uniformly good policy implies the sub-linear growth of the
system regret and achieves the maximum long-term average reward
$\Sigma_{j=1}^M\mu(\theta_{\sigma(j)})$ which is the same as in the
case with perfect knowledge of $\Theta$.

\subsection{The Logarithmic Order and the Optimal Policy}
We present in the theorem below the result established in [1,2] on
the logarithmic order as well as the leading constant of the minimum
regret growth rate of the single-player MAB.

{\em Theorem~[1,2]:} Under the regularity conditions (conditions
C1-C4 in Appendix~A), we have, for any uniformly good policy $\pi$,
\begin{eqnarray}\label{eqn:cenlow}
\liminf_{T\rightarrow\infty}\frac{R_T^{\pi}(\Theta)}{\log
T}\ge\Sigma_{j:~\mu(\theta_j)<\mu(\theta_{\sigma(M)})}\frac{\mu(\theta_{\sigma(M)})-\mu(\theta_{j})}
{I(\theta_{j},\theta_{\sigma(M)})}.
\end{eqnarray}

\vspace{1em}

Lai and Robbins also constructed a policy that achieves the lower
bound on the regret growth rate given in~\eqref{eqn:cenlow} under
single play ($M=1$)~\cite{Lai&Robbins85AAM} (which was extended by
Anantharam \etal to $M>1$ in~\cite{Anantharam&etal87TAC}). Under
this policy, two statistics of the past observations are maintained
for each arm. Referred to as the point estimate, the first statistic
is an estimate of the mean $\mu(\theta_i)$ given by a function $h_k$
of the past observations on this arm ($k$ denotes the total number
of observations). The second statistic is the so-called confidence
upper bound which represents the \emph{potential} of an arm: the
less frequently an arm is played, the less confident we are about
the point estimate, and the higher the potential of this arm. The
confidence upper bound, denoted by $g_{t,k}$, thus depends on not
only the number $k$ of observations on the arm but also the current
time $t$ in order to measure how frequently this arm has been
played.

\begin{figure}[htbp]
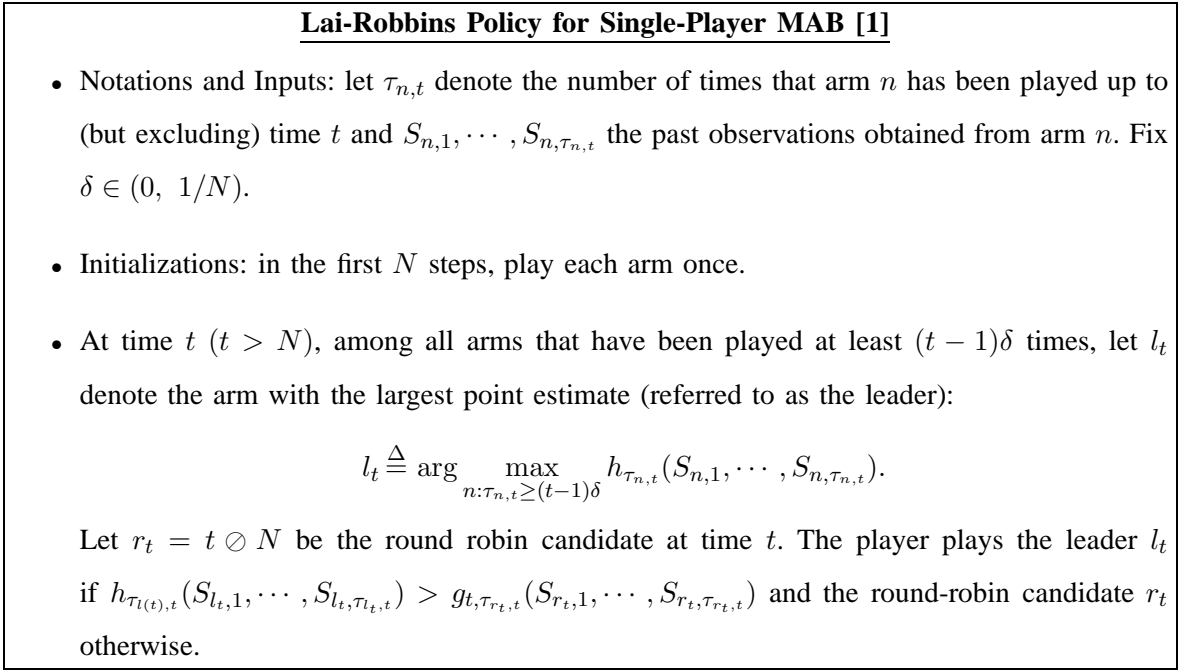

\begin{center}
\noindent\fbox{
\parbox{6in}
{ \centerline{\underline{{\bf Lai-Robbins Policy for Single-Player
MAB~\cite{Lai&Robbins85AAM}}}} {
\begin{itemize}
\item Notations and Inputs:  let $\tau_{n,t}$ denote the number of
times that arm $n$ has been played up to (but excluding) time $t$
and $S_{n,1},\cdots,S_{n,\tau_{n,t}}$ the past observations obtained
from
arm $n$. Fix $\delta\in(0,~1/N)$.\\[-1em]
\item Initializations: in the first $N$ steps, play each arm
once.\\[-1em]
\item At time $t~(t>N)$, among all arms that have been played at least
$(t-1)\delta$ times, let $l_t$ denote the arm with the largest point
estimate (referred to as the
leader):\[l_t\defeq\arg\max_{n:\tau_{n,t}\ge(t-1)\delta}h_{\tau_{n,t}}(S_{n,1},\cdots,S_{n,\tau_{n,t}}).\]
Let $r_t=t\oslash N$ be the round robin candidate at time $t$. The
player plays the leader $l_t$ if
$h_{\tau_{l(t),t}}(S_{l_t,1},\cdots,S_{l_t,\tau_{l_t,t}})>
g_{t,\tau_{r_t,t}}(S_{r_t,1},\cdots,S_{r_t,\tau_{r_t,t}})$ and the
round-robin candidate $r_t$ otherwise.
\end{itemize}}}}
\caption{Lai-Robbins policy for single-player
MAB~\cite{Lai&Robbins85AAM}.}\label{fig:LaiRobbinPolicy}
\end{center}
\end{figure}

Based on these two statistics, Lai-Robbins policy operates as
follows. At each time $t$, among all ``well-sampled'' arms, the one
with the largest point estimate is selected as the leader denoted as
$l_t$. The player then chooses between the leader $l_t$ and a
round-robin candidate $r_t=t\oslash N$ to play. The leader $l_t$ is
played if and only if its point estimate exceeds the confidence
upper bound of the round-robin candidate $r_t$. A detailed
implementation of this policy is given in
Fig.~\ref{fig:LaiRobbinPolicy}.

Lai and Robbins~\cite{Lai&Robbins85AAM} have shown that for point
estimates $h_k$ and confidence upper bounds $g_{t,k}$ satisfying
certain conditions (condition C5 in Appendix~A), the above policy is
optimal, \ie it achieves the minimum regret growth rate given
in~\eqref{eqn:cenlow}. For Gaussian, Bernoulli, and Poisson reward
models, $h_k$ and $g_{t,k}$ satisfying condition C5 are given
below~\cite{Lai&Robbins85AAM}.
\begin{eqnarray}\label{eqn:pt}
h_{k}&=&\frac{\sum_{j=1}^kS_{n,j}}{k},\\\label{eqn:cub}
g_{t,k}&=&\inf\{\lambda:~\lambda\ge
h_{k}~\mbox{and}~I(h_k,\lambda)\ge\frac{\log(t-1)}{k}\}.
\end{eqnarray}While $g_{t,k}$ given in~\eqref{eqn:cub} is not in closed-form, it does
not affect the implementation of the policy. Specifically, the
comparison between the point estimate of the leader $l_t$ and the
confidence upper bound of the round-robin candidate $r_t$ (see
Fig.~\ref{fig:LaiRobbinPolicy}) is shown to be equivalent to the
following two conditions~\cite{Lai&Robbins85AAM}:
{\begin{eqnarray}\nn
&&h_{\tau_{r_t,t}}(S_{r_t,1},\cdots,S_{r_t,\tau_{r_t,t}})<h_{\tau_{l_t,t}}(S_{l_t,1},\cdots,S_{l_t,\tau_{l_t,t}})\\\nn
&&I(h_{\tau_{r_t,t}},h_{\tau_{l_t,t}})>\frac{\log(t-1)}{\tau_{r_t,t}}.
\end{eqnarray}}
Consequently, we only require the point estimate $h_{k}$ to
implement the policy.

\subsection{Order-Optimal Index Policies}
Since Lai and Robbins's seminal work, researchers have developed
several index-type policies that are simpler than Lai-Robbins policy
by using a single sample mean based statistic [8,9]. Specifically,
under such an index policy, each arm is assigned with an index that
is a function of its sample mean and the arm with the greatest index
will be played in each slot. To obtain an initial index, each arm is
played once in the first $N$ slots.

The indexes proposed in~\cite{Agrawal95AAP} for Gaussian, Bernoulli,
Poisson and exponential distributed reward models are given
in~\eqref{eqn:index}. Except for Gaussian distribution, this index
policy only achieves the optimal logarithmic order of the regret
growth rate but not the best leading constant\footnote{The optimal
indexes for Bernoulli, Poisson and exponential distributions are
also developed in~\cite{Agrawal95AAP}. However, the indexes are not
given in closed-form and are difficult to implement.}.
\begin{eqnarray}\label{eqn:index}
I(x)=\left\{\begin{array}{ll}
x+(2\frac{\log(t-1)+2\log\log(t-1)}{\tau_{n,t}})^{1/2},
& \mbox{Gaussian}\\
x+\min\{(2\frac{\log(t-1)+2\log\log(t-1)}{\tau_{n,t}})^{1/2}/2,~1\}, & \mbox{Bernoulli}\\
x+\min\{(2a\frac{\log(t-1)+2\log\log(t-1)}{\tau_{n,t}})^{1/2}/2,~a\},& \mbox{Poisson}\\
x+b\min\{(2\frac{\log(t-1)+2\log\log(t-1)}{\tau_{n,t}})^{1/2},~1\},
& \mbox{Exponential}
\end{array}\right.,
\end{eqnarray}
where $I(\cdot)$ is the index function, $x$ is the sample mean of
the arm, and $a,b$ are upper bounds of all possible values of
$\theta$ in the Poisson and exponential reward models, respectively.

Based on~\cite{Agrawal95AAP}, a simpler order-optimal sample mean
based index policy was established in~\cite{Auer&etal02ML} for
reward distributions with a known finite support:
\begin{eqnarray}\label{eqn:auerIndex}
I(x)=x+(2\frac{\log(t-1)}{\tau_{n,t}})^{1/2}.
\end{eqnarray}

\section{Decentralized MAB: Problem Formulation}\label{sec:problemformulation}
In this section, we formulate the decentralized MAB with $M$
distributed players. In addition to conditions C1-C5 required by the
centralized MAB, we assume that the $M$ best arms have distinct
nonnegative means\footnote{This condition can be relaxed to the case
that a tie occurs at the $M$th largest mean.}.

In the decentralized setting, players may collide and may not
receive reward that the arm can potentially offer. We thus refer to
$S_i(t)$ as the \emph{state} of arm $i$ at time $t$ (for example,
the busy/idle state of a communication channel in the context of
cognitive radio). At time $t$, player $i~(1\le i\le M)$ chooses an
action $a_i(t)\in\{1,\cdots,N\}$ that specifies the arm to play and
observes its state $S_{a_i(t)}(t)$. The action $a_i(t)$ is based on
the player's \emph{local} observation and decision history. Note
that the local observation history of each player also includes the
collision history. As shown in Sec.~\ref{subsec:ElimPre}, the
observation of a collision is used to adjust the local offset of a
player's time sharing schedule of the $M$ best arms to avoid
excessive future collisions.

We define a local policy $\pi_i$ for player $i$ as a sequence of
functions $\pi_i=\{\pi_i(t)\}_{t\ge1}$, where $\pi_i(t)$ maps player
$i$'s past observations and decisions to action $a_i(t)$ at time
$t$. The decentralized policy $\pi$ is thus given by the
concatenation of the local polices for all players:
{\[\pi\defeq[\pi_1,\cdots,\pi_M].\]} Define reward $Y(t)$ as the
total reward accrued by all players at time $t$, which depends on
the system collision model as given below.

{\em Collision model 1:} When multiple players choose the same arm
to play, they share the reward in an arbitrary way. Since how
players share the reward has no effect on the total system regret,
without loss of generality, we assume that only one of the colliding
players obtains a random reward given by the current state of the
arm. Under this model, we have
{\[Y(t)=\Sigma_{j=1}^N\mathbb{I}_j(t)S_j(t),\]}where
$\mathbb{I}_j(t)$ is the indicator function that equals to $1$ if
arm $j$ is played by \emph{at least} one player, and $0$ otherwise.

{\em Collision model 2:} When multiple players choose the same arm
to play, no one obtains reward. Under this model, we have
{\[Y(t)=\Sigma_{j=1}^N\hat{\mathbb{I}}_j(t)S_j(t),\]}where
$\hat{\mathbb{I}}_j(t)$ is the indicator function that equals to $1$
if arm $j$ is played by \emph{exactly} one player, and $0$
otherwise.

The system regret of policy $\pi$ is thus given by
{\[R_T^{\pi}(\Theta)=T\Sigma_{j=1}^M\mu(\theta_{\sigma(j)})-\mathbb{E}_{\pi}[\Sigma_{t=1}^TY(t)].\]}
Note that the system regret in the decentralized MAB is defined with
respect to the same best possible reward
$T\Sigma_{j=1}^M\mu(\theta_{\sigma(j)})$ as in its centralized
counterpart. The objective is to minimize the rate at which
$R_T(\Theta)$ grows with time $T$ under any parameter set $\Theta$
by choosing an optimal decentralized policy. Similarly, we say a
decentralized policy is {\em uniformly good} if for any parameter
set $\Theta$, we have $R_T(\Theta)=o(T^b)$ for any $b>0$. To address
the optimal order of the regret, it is sufficient to focus on
uniformly good decentralized polices provided that such policies
exist.

We point out that all results developed in this work apply to a more
general observation model. Specifically, the arm state observed by
different players can be drawn from different distributions, as long
as players have the common set of the $M$ best arms and each of the
$M$ best arms has the same mean across players. This relaxation in
the observation model is particularly important in the application
of opportunistic transmission in fading channels, where different
users experience different fading environments in the same channel.

\section{The Optimal Order of The System Regret}\label{sec:orderoptimal}
In this section, we show that the optimal order of the system regret
growth rate of the decentralized MAB is logarithmic, the same as its
centralized counterpart as given in Sec.~\ref{sec:classic}.

\begin{theorem}\label{thm:optorder}
Under both collision models, the optimal order of the system regret
growth rate of the decentralized MAB is logarithmic, \ie for an
optimal decentralized policy $\pi^*$, we have
{\begin{eqnarray}\label{eqn:optorder}
L(\Theta)\le\liminf_{T\rightarrow\infty}\frac{R_T^{\pi^*}(\Theta)}{\log
T}\le\limsup_{T\rightarrow\infty}\frac{R_T^{\pi^*}(\Theta)}{\log
T}\le U(\Theta)
\end{eqnarray}}
for some constants $L(\Theta)$ and $U(\Theta)$ that depend on
$\Theta$.
\end{theorem}
\begin{proof}
The proof consists of two parts. First, we prove that the lower
bound for the centralized MAB given in~\eqref{eqn:cenlow} is also a
lower bound for the decentralzied MAB. Second, we construct a
decentralized policy (see Sec.~\ref{sec:denpolicy}) that achieves
the logarithmic order of the regret growth rate.

It appears to be intuitive that the lower bound for the centralized
MAB provides a lower bound for the decentralized MAB. This, however,
may not hold when some of the $M$ best arms have negative means
(modeling the punishment for playing certain arms). The reason is
that the centralized MAB considered in~\cite{Anantharam&etal87TAC}
requires exactly $M$ arms to be played at each time, while in the
decentralized setting, fewer than $M$ arms may be played at each
time when players choose the same arm. To obtain an lower bound for
the decentralized MAB, we need to consider a centralized MAB where
the player has the freedom of playing \emph{up to} $M$ arms at each
time. When the player knows that \emph{all} arms have nonnegative
means, it is straightforward to see that the two conditions of
playing exactly $M$ arms and playing no more than $M$ arms are
equivalent. Without such knowledge, however, this statement needs to
be proven.

\begin{lemma}\label{lemma:equiv}
Under the condition that the $M$ best arms have nonnegative means,
the centralized MAB that requires exactly $M$ arms are played at
each time is equivalent to the one that requires at most $M$ arms
are played at each time.
\end{lemma}
\emph{Proof:} See Appendix~B.
\end{proof}

\section{An Order-Optimal Decentralized Policy}\label{sec:denpolicy}
In this section, we construct a decentralized policy $\pi^*_F$ that
achieves the optimal logarithmic order of the system regret growth
rate under the fairness constraint.

\subsection{Basic Structure of the Decentralized TDFS Policy
$\pi^*_F$}\label{subsec:DecenStructure} The basic structure of the
proposed policy $\pi^*_F$ is a time division structure at each
player for selecting the $M$ best arms. For the ease of the
presentation, we first assume that there is a pre-agrement among
players so that they use different phases (offsets) in their time
division schedule. For example, the offset in each player's time
division schedule can be predetermined based on the player's ID. In
Sec.~\ref{subsec:ElimPre}, we show that this pre-agreement can be
eliminated while maintaining the order-optimality and fairness of
the TDFS policy, which leads to a complete decentralization among
players.

Consider, for example, the case of $M=2$. The time sequence is
divided into two disjoint subsequences, where the first subsequence
consists of all odd slots and the second one consists of all even
slots. The pre-agreement is such that player $1$ targets at the best
arm during the first subsequence and the second best arm during the
second subsequence, and player $2$ does the opposite.

\begin{figure}[h]
\centerline{
\begin{psfrags}
\psfrag{a1}[c]{\footnotesize$~\hat{\alpha}(1)=1$}\psfrag{a2}[c]{\footnotesize$~~\hat{\alpha}(1)=2$}
\psfrag{a3}[c]{\footnotesize$~~\hat{\alpha}(1)=3$}\psfrag{t=1}[c]{\small$t=1$}\psfrag{2}[c]{\small$2$}
\psfrag{3}[c]{\small$3$}\psfrag{4}[c]{\small$4$}\psfrag{5}[c]{\small$5$}
\psfrag{6}[c]{\small$6$}\psfrag{7}[c]{\small$7$}\psfrag{8}[c]{\small$8$}\psfrag{9}[c]{\small$9$}
\psfrag{10}[c]{\small$10$}\psfrag{11}[c]{\small$11$}\psfrag{12}[c]{\small$12$}
\psfrag{s1}[l]{The first subsequence (targeting at the best arm)}
\psfrag{m1}[l]{The first mini-sequence (targeting at the second best
arm when arm~$1$ is considered as the best)}\psfrag{m2}[l]{The
second mini-sequence (targeting at the second best arm when arm~$2$
is considered as the best)}\psfrag{m3}[l]{The third mini-sequence
(targeting at the second best arm when arm~$3$ is considered as the
best)}\scalefig{1}\epsfbox{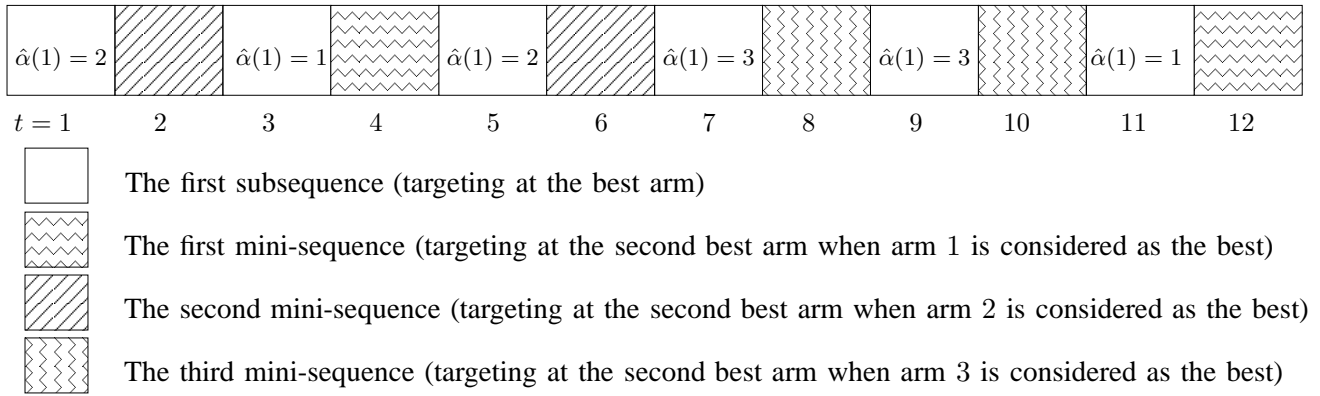}
\end{psfrags}}
\caption{The structure of player 1's local policy under
$\pi^*_F$~($M=2$,~$N=3$. Let $\hat{\alpha}(1)$ denote the arm
considered as the best by player 1 in the first subsequence. In this
example, player 1 divides the second subsequence (\ie all the even
slots) into three mini-sequences, each associated with a subset of
$2$ arms after removing the arm considered to be the best in the
first subsequence.).} \label{fig:sequence}
\end{figure}

Without loss of generality, consider player~1. In the first
subsequence, player 1 applies a single-player policy, say
Lai-Robbins policy, to efficiently learn and select the best arm. In
the second subsequence, the second best arm is learned and
identified by removing the arm that is considered to be the best and
applying the Lai-Robbins policy to identify the best arm in the
remaining $N-1$ arms (which is the second best among all $N$ arms).
Note that since the best arm cannot be perfectly identified, which
arm is considered as the best in an odd slot $t$ (in the first
subsequence) is a random variable determined by the realization of
past observations. We thus partition the second subsequence into
multiple mini-sequences depending on which arm was considered to be
the best in the preceding odd slot and thus should be removed from
consideration when identifying the second best. Specifically, as
illustrated in Fig.~\ref{fig:sequence}, the second subsequence is
divided into $N$ disjoint mini-sequences, where the $i$th
mini-sequence consists of all slots that follow a slot in which arm
$i$ was played (\ie arm $i$ was considered the best arm in the
preceding slot that belongs to the first subsequence). In the $i$th
mini-sequence, the player applies Lai-Robbins policy to arms
$\{1,\cdots, i-1, i+1,\cdots, N\}$ after removing arm $i$.

In summary, the local policy for each player consists of $N+1$
parallel Lai-Robbins procedures: one is applied in the subsequence
that targets at the best arm and the rest $N$ are applied in the $N$
mini-sequences that target at the second best arm. These parallel
procedures, however, are coupled through the common observation
history since in each slot, regardless of which subsequence or
mini-sequence it belongs to, all the past observations are used in
the decision making. We point out that making each mini-sequence use
only its own observations is sufficient for the order optimality of
the TDFS policy and simplifies the optimality analysis. However, we
expect that using all available observations leads to a better
constant as demonstrated by simulation results in
Sec.~\ref{sec:simu}.

The basic structure of the TDFS policy $\pi^*_F$ is the same for the
general case of $M>2$. Specifically, the time sequence is divided
into $M$ subsequences, in which each player targets at the $M$ best
arms in a round-robin fashion with a different offset. Suppose that
player 1 has offset $0$, \ie it targets at the $k$th ($1\le k\le M$)
best arm in the $k$th subsequence. To player~1, the $k$th~$(k>1)$
subsequence is then divided into $C^N_{k-1}$ mini-sequences, each
associated with a subset of $N-k+1$ arms after removing those $k-1$
arms that are considered to have a higher rank than $k$. In each of
the $C^N_{k-1}$ mini-sequences, Lai-Robbins single-player policy is
applied to the subset of $N-k+1$ arms associated with this
mini-sequence. Note that while the subsequences are deterministic,
each mini-sequence is random. Specifically, for a given slot $t$ in
the $k$th subsequence, the mini-sequence which it belongs to is
determined by the specific actions of the player in the previous
$k-1$ slots. For example, if arms $1,\cdots,k-1$ are played in the
previous $k-1$ slots, then slot $t$ belongs to the mini-sequence
associated with the arm set of $\{k,\cdots,N\}$. A detailed
implementation of the proposed policy for a general $M$ is given in
Fig.~\ref{fig:decenpolicy}.

Note that the basic structure of the TDFS policy $\pi^*_F$ is
general. It can be used with any single-player policy to achieve an
efficient and fair sharing of the $M$ best arms among $M$
distributed players. Furthermore, players are not restricted to
using the same single-player policy. Details are given in
Sec.~\ref{subsec:orderoptimal}.

\begin{figure}[htbp]
\begin{center}
\noindent\fbox{
\parbox{6in}
{ \centerline{\underline{{\bf The Decentralized TDFS Policy
$\pi^*_F$}}} {Without loss of generality, consider player $i$.
\begin{itemize}
\item Notations and Inputs: In addition to the notations and inputs of Lai-Robbins
single-player policy (see Fig.~\ref{fig:LaiRobbinPolicy}), let
$m_k(t)$ denote the number of slots in the $k$th subsequence up to
(and including) $t$, and let $m_{\Nc}(t)$ denote the number of slots
in the
mini-sequence associated with arm set $\Nc$ up to (and including) $t$.\\[-1em]
\item At time $t$, player $i$ does the following.
\begin{enumerate}
\item[1.] If $t$ belongs to the $i$th subsequence (\ie $t\oslash M=i$), player $i$ targets at the best arm
by carrying out the following procedure.\\
If $m_i(t)\le N$, play arm $m_i(t)$. Otherwise, the player chooses
between a leader $l_t$ and a round-robin candidate
$r_t=m_i(t)\oslash N$, where the leader $l_t$ is the arm with the
largest point estimate among all arms that have been played for at
least $(m_i(t)-1)\delta$ times. The player plays the leader $l_t$ if
its point estimate is larger than the confidence upper bound of
$r_t$. Otherwise, the player plays the round-robin candidate $r_t$.\\[-1em]
\item[2.] If $t$ belongs to the $k$th ($k\neq i$) subsequence (\ie $t\oslash
M=k$), the player targets at the $j$th best arm where
$j=(k-i+M+1)\oslash M$ by carrying out the following procedure.\\
If $m_k(t)\le N$, play arm $m_k(t)$. Otherwise, the player targets
at the $j$th best arm. Let $\Ac_t$ denote the set of $j-1$ arms
played in the previous $j-1$ slots. Slot $t$ thus belongs to the
mini-sequence associated with the subset
$\Nc_t=\{1,\cdots,N\}\backslash\Ac_t$ of arms. The player chooses
between a leader and a round-robin candidate defined within $\Nc_t$.
Specifically, among all arms that have been played for at least
$(m_{\Nc_t}(t)-1)\delta$ times, let $l_t$ denote the leader. Let
$r_t=m_{\Nc_t}(t)\oslash (N-j+1)$ be the round-robin candidate
where, for simplicity, we have assumed that arms in $\Nc_t$ are
indexed by $1,\cdots,N-j+1$. The player plays the leader $l_t$ if
its point estimate is larger than the confidence upper bound of
$r_t$. Otherwise, the player plays the round-robin candidate $r_t$.
\end{enumerate}
\end{itemize}} }} \caption{The decentralized TDFS policy
$\pi^*_F$.}\label{fig:decenpolicy}
\end{center}
\end{figure}

\subsection{Order-Optimality under Fairness
Constraint}\label{subsec:orderoptimal} Compared to the single-player
counterpart given in~[1,2], the difficulties in establishing the
logarithmic regret growth rate of $\pi^*_F$ are twofold. First,
compared to the centralized problem where a single player observes
$M$ different arms simultaneously, each player here can only observe
one arm at each time. Each player is thus learning the entire rank
of the $M$ best arms with fewer observations. Furthermore, since the
rank of any arm can never be perfectly identified, the mistakes in
identifying the $i$th ($1\le i< M$) best arm will propagate into the
learning process for identifying the $(i+1)$th, up to the $M$th best
arm. Second, without centralized scheduling, collisions are bound to
happen even with pre-agreed offset on sharing the $M$ best arms
since players do not always agree on the rank of the arms. Such
issues need to be carefully dealt with in establishing the
logarithmic order of the regret growth rate.

\begin{theorem}\label{thm:worstcase}
Under the TDFS policy $\pi^*_F$, we have, for some constant
$C(\Theta)$, {\begin{eqnarray}\label{eqn:up}
\limsup_{T\rightarrow\infty}\frac{R^{\pi^*_F}_T(\Theta)}{\log T}\le
C(\Theta).
\end{eqnarray}}Let
\[
x_k\defeq\Sigma_{i=1}^k\Sigma_{j:\mu(\theta_j)<\mu(\theta_{\sigma(i)})}
\frac{1}{I(\theta_{j},\theta_{\sigma(i)})}.
\]
Under collision model 1, {\begin{eqnarray}\nn
C(\Theta)=M(\Sigma_{i=1}^Mx_k\mu(\theta_{\sigma(i)})-\Sigma_{n:\mu(\theta_n)<\mu(\theta_{\sigma(M)})}
\frac{\mu(\theta_{n})}{I(\theta_{n},\theta_{\sigma(M)})}).\nn
\end{eqnarray}}Under collision model 2,
{\begin{eqnarray}\nn C(\Theta)=
M(\Sigma_{i=1}^M\Sigma_{k=1}^Mx_k\mu(\theta_{\sigma(i)})-\Sigma_{n:\mu(\theta_n)<\mu(\theta_{\sigma(M)})}
\mu(\theta_{n})\max\{\frac{1}{I(\theta_{n},\theta_{\sigma(M)})}-\Sigma_{i=1}^{M-1}
\frac{1}{I(\theta_{n},\theta_{\sigma(i)})},0\}).
\end{eqnarray}}
\end{theorem}
\begin{proof}
Note that the system regret is given by the sum of the regrets in
the $M$ subsequences. By symmetry in the structure of $\pi^*_F$, all
subsequences experience the same regret. Consequently, the system
regret is equal to $M$ times the regret in each subsequence. We thus
focus on a particular subsequence and show that the regret growth
rate in this subsequence is at most logarithmic. Without loss of
generality, consider the first subsequence in which the $i$th player
targets at the $i$th ($1\le i\le M$) best arm.

To show the logarithmic order of the regret growth rate, The basic
approach is to show that the number of slots in which the $i$th
player does not play the $i$th best arm is at most logarithmic with
time. This is done by establishing in Lemma~\ref{lemma:lbbesttime} a
lower bound on the number of slots in which the $i$th player plays
the $i$th best arm in the first subsequence.

To show Lemma~\ref{lemma:lbbesttime}, we first establish
Lemma~\ref{lemma:upbadtime} by focusing on the \emph{dominant}
mini-sequence of the first subsequence. Without loss of generality,
consider player~$i$. Define the dominant mini-sequence as the one
associated with arm set $\{\sigma(i),\cdots,\sigma(N)\}$ (\ie the
$i-1$ best arms are correctly identified and removed in this
mini-sequence). Lemma~\ref{lemma:upbadtime} shows that, in the
dominant mini-sequence, the number of slots in which player $i$ does
not play the $i$th best arm is at most logarithmic with time.

\begin{lemma}\label{lemma:upbadtime}
Let $\tau_{n,T}^D$ denote the number of slots up to time $T$ in
which arm $n$ is played in the dominant mini-sequence associated
with arm set $\{\sigma(i),\cdots,\sigma(N)\}$. Then, for any arm $n$
with $\mu(\theta_n)<\mu(\theta_{\sigma(i)})$, we have,
\begin{eqnarray}
\limsup_{T\rightarrow\infty}\frac{\mathbb{E}[\tau_{n,T}^D]}{\log
T}\le\frac{1}{I(\theta_n,\theta_{\sigma(i)})}.
\end{eqnarray}
\end{lemma}
\begin{proof}
Note that this lemma is an extension of Theorem 3
in~\cite{Lai&Robbins85AAM}. The proof of this lemma is more
complicated since the mini-sequence is random and the decisions made
in this mini-sequence depend on all past observations (no matter to
which mini-sequence they belong). See Appendix~C for details.
\end{proof}

Next, we establish Lemma~\ref{lemma:lbbesttime}. The basic approach
is to show that the length of the dominant mini-sequence dominates
the lengths of all other mini-sequences in the first subsequences.
Specifically, we show that the number of slots that do not belong to
the dominant mini-sequence is at most logarithmic with time. This,
together with Lemma~\ref{lemma:upbadtime} that characterizes the
dominant mini-sequence, leads to Lemma~\ref{lemma:lbbesttime} below.

\begin{lemma}\label{lemma:lbbesttime}
Let $\hat{\tau}_{\sigma(i),T}$ denote the number of slots in which
player $i$ plays the $i$th best arm in the first subsequence up to
time $T$, we have
\begin{eqnarray}\label{eqn:lb1}
\limsup_{T\rightarrow\infty}\frac{T/M-\mathbb{E}[\hat{\tau}_{\sigma(i),T}]}{\log
T}\le x_i,
\end{eqnarray}where
\begin{eqnarray}
x_i=\Sigma_{j=1}^i\Sigma_{k:\mu(\theta_k)<\mu(\theta_{\sigma(j)})}\frac{1}{I(\theta_{\sigma(k)},\theta_{\sigma(j)})}.
\end{eqnarray}
\end{lemma}
\begin{proof} The proof is based on an induction argument on $i$.
See Appendix~D for details.
\end{proof}

From Lemma~\ref{lemma:lbbesttime}, for all $1\le i\le M$, the number
of slots that the $i$th best arm is not played by player $i$ is at
most logarithmic with time. Consequently, for all $j\neq i$, the
number of slots that player $j$ plays the $i$th best arm is also at
most logarithmic with time, \ie the number of collisions on the
$i$th best arm is at most logarithmic with time. Since a reward loss
on the $i$th $(1\le i\le M)$ best arm can only happen when it is not
played or a collision happens, the reward loss on the $i$th best is
at most logarithmic order of time, leading to the logarithmic order
of the regret growth rate.

To establish the upper bound $C(\Theta)$ on the constant of the
logarithmic order of the system regret growth rate, we consider the
worst-case collisions on each arm. See Appendix~E for details.
\end{proof}

From Theorem~\ref{thm:optorder} and Theorem~\ref{thm:worstcase}, the
decentralized policy is order-optimal. Furthermore, as given in
Theorem~\ref{thm:fair} below, the decentralized policy ensures
fairness among players under a fair collision model that offers all
colliding players the same \emph{expected} reward. For example, if
only one colliding player can receive reward, then all players have
equal probability to be lucky.
\begin{theorem}\label{thm:fair}
Define the local regret for player $i$ under $\pi^*_F$ as
{\[R^{\pi^*_F}_{T,i}(\Theta)\defeq
\frac{1}{M}T\Sigma_{j=1}^M\mu(\theta_{\sigma(j)})-\mathbb{E}_{\pi^*_F}[\Sigma_{t=1}^TY_i(t)],
\]}where $Y_i(t)$ is the immediate reward obtained by player $i$ at time
$t$. Under a fair collision model, we have,
{\begin{eqnarray}\label{eqn:fair}
\limsup_{T\rightarrow\infty}\frac{R^{\pi^*_F}_{T,i}(\Theta)}{\log
T}=\frac{1}{M}\limsup_{T\rightarrow\infty}\frac{R^{\pi^*_F}_T(\Theta)}{\log
T}.
\end{eqnarray}}
\end{theorem}

Theorem~\ref{thm:fair} follows directly from the symmetry among
players under $\pi^*_F$. It shows that each player achieves the same
time average reward
$\frac{1}{M}T\Sigma_{j=1}^M\mu(\theta_{\sigma(j)})$ at the same
rate. We point out that without knowing the reward rate that each
arm can offer, ensuring fairness requires that each player identify
the entire set of the $M$ best arms and share each of these $M$ arms
evenly with other players. As a consequence, each player needs to
learn which of the $C^N_M$ possibilities is the correct choice. This
is in stark difference from polices that make each player target at
a single arm with a specific rank (for example, the $i$th player
targets solely at the $i$th best arm). In this case, each player
only needs to distinguish one arm (with a specific rank) from the
rest. The uncertainty facing each player, consequently, the amount
of learning required, is reduced from $C^N_M$ to $N$. Unfortunately,
fairness among players is lost.

As mentioned in Sec.~\ref{subsec:DecenStructure}, the proposed
policy $\pi^*_F$ can be used with any single-player policy, which
can also be different for different players. More importantly, the
order optimality of the TDFS policy is preserved as long as each
player's single-player policy achieves the optimal logarithmic order
in the single player setting. This statement can be proven along the
same line as Theorem~\ref{thm:worstcase} by establishing results
similar to Lemma~\ref{lemma:upbadtime} and
Lemma~\ref{lemma:lbbesttime}.

\subsection{Eliminating the Pre-Agreement}\label{subsec:ElimPre}

In this subsection, we show that a pre-agreement among players on
the offsets of the time division schedule can be eliminated in the
TDFS policy $\pi_F^*$ while maintaining its order-optimality and
fairness.

Specifically, when a player joins the system, it randomly generates
a local offset uniformly drawn from $\{0,1,\cdots,M-1\}$ and plays
one round of the $M$ arms considered to be the best. For example, if
the random offset is $1$, the player targets at the second, the
third, $\cdots$, the $M$th, and then the best arms in the $M$
subsequent slots, respectively. If no collision is observed in this
round of $M$ plays, the player keeps the same offset; otherwise the
player randomly generates a new offset. Over the sequence of slots
where the same offset is used, the player implements the local
policy of $\pi_F^*$ (given in Fig.~\ref{fig:decenpolicy}) with this
offset. To summarize, each player implements $M$ parallel local
procedures of $\pi_F^*$ corresponding to $M$ different offsets.
These parallel procedures are coupled through the observation
history, \ie each player uses its entire local observation history
in learning no matter which offset is being used.

Note that players can join the system at different time. We also
allow each player to leave the system for an arbitrary finite number
of slots.
\begin{theorem}\label{thm:optElimPre} The decentralized TDFS
policy $\pi_F^*$ without pre-agreement is order-optimal and fair.
\end{theorem}
\begin{proof} See Appendix~F.
\end{proof}

\section{A Lower Bound for a Class of Decentralized
Polices}\label{sec:lowerbound} In this section, we establish a lower
bound on the growth rate of the system regret for a general class of
decentralized policies, to which the proposed policy $\pi^*_F$
belongs. This lower bound provides a tighter performance benchmark
compared to the one defined by the centralized MAB. The definition
of this class of decentralized polices is given below.

\begin{definition}{\em Time Division Selection Policies}
The class of time division selection (TDS) policies consists of all
decentralized polices $\pi=[\pi_1,\cdots,\pi_M]$ that satisfy the
following property: under local policy $\pi_i$, there exists $0\le
a_{ij}\le 1~(1\le i,j\le M,~\Sigma_{j=1}^Ma_{ij}=1~\forall~i)$
independent of the parameter set $\Theta$ such that the expected
number of times that player $i$ plays the $j$th $(1\le j\le M)$ best
arm up to time $T$ is $a_{ij}T-o(T^b)$ for all $b>0$.
\end{definition}

A policy in the TDS class essentially allows a player to efficiently
select each of the $M$ best arms according to a fixed time portion
that does not depend on the parameter set $\Theta$. It is easy to
see that the TDFS policy $\pi_F^*$ (with or without pre-agreement)
belongs to the TDS class with $a_{ij}=1/M$ for all $1\le i,~j\le M$.

\begin{theorem}\label{thm:lowerboundTDS}
For any uniformly good decentralized policy $\pi$ in the TDS class,
we have {\begin{eqnarray}\label{eqn:delow}
\liminf_{T\rightarrow\infty}\frac{R_T^{\pi}(\Theta)}{\log
T}\ge\Sigma_{i=1}^M\Sigma_{j:\mu(\theta_j)<\mu(\theta_{\sigma(M)})}\frac{\mu(\theta_{\sigma(M)})-\mu(\theta_{j})}
{I(\theta_{j},\theta_{\sigma(i)})}.
\end{eqnarray}}
\end{theorem}
\begin{proof}
The basic approach is to establish a lower bound on the number of
slots in which each player plays an arm that does not belong to the
$M$ best arms. By considering the best case that they do not
collide, we arrive at the lower bound on the regret growth rate
given in~\eqref{eqn:delow}. The proof is based on the following
lemma, which generalizes Theorem 2 in~\cite{Lai&Robbins85AAM}. To
simplify the presentation, we assume that the means of all arms are
distinct. However, Theorem~\ref{thm:lowerboundTDS} and
Lemma~\ref{lemma:lowboundgen} apply without this assumption.
\begin{lemma}\label{lemma:lowboundgen} Consider a local policy $\pi_i$. If for any
parameter set $\Theta$ and $b>0$, there exists a
$\Theta$-independent positive increasing function $v(T)$ satisfying
$v(T)\rightarrow\infty$ as $T\rightarrow\infty$ such that
\begin{eqnarray}\label{eqn:lowtime}
\mathbb{E}[v(T)-\tau_{\sigma(j),T}]\le o((v(T))^b),
\end{eqnarray}
then we have, $\forall~k>j$,
\begin{eqnarray}\label{eqn:l}
\liminf_{T\rightarrow\infty}\frac{\mathbb{E}[\tau_{\sigma(k),T}]}{\log
v(T)}\ge \frac{1}{I(\theta_{\sigma(k)},\theta_{\sigma(j)})}.
\end{eqnarray}
\end{lemma}

Note that Lemma~\ref{lemma:lowboundgen} is more general than Theorem
2 in~\cite{Lai&Robbins85AAM} that assumes $v(T)=T$ and $j=1$. The
proof of Lemma~\ref{lemma:lowboundgen} is given in Appendix~G.

Consider a uniformly good decentralized policy $\pi$ in the TDS
class. There exists a player, say player $u_1$, that plays the best
arm for at least $T/M-o(T^b)$ times. Since at these times, $u_1$
cannot play other arms, there must exist another player, say player
$u_2$, that plays the second best arm for at least
$T/(M(M-1))-o(T^b)$ times. It thus follows that there exist $M$
different players $u_1,\cdots,u_M$ such that under any parameter set
$\Theta$, the expected time player $u_i~(1\le i\le M)$ plays the
$i$th best arm is at least $T/\prod_{j=1}^i(M-j+1)-o(T^b)$ times.
Based on Lemma~\ref{lemma:lowboundgen}, for any arm $j$ with
$\mu(\theta_j)<\mu(\theta_{\sigma(M)})$, the expected time that
player $u_i$ plays arm $j$ is at least
$1/I(\theta_j;\theta_{\sigma(i)})\log
(T/\prod_{j=1}^i(M-i+1))-o(\log T)$. By considering the best case
that players do not collide, we arrive at~\eqref{eqn:delow}.

\end{proof}

\section{Simulation Examples}\label{sec:simu}
In this section, we study the performance (\ie the leading constant
of the logarithmic order) of the decentralized TDFS policy in
different applications through simulation examples.

\subsection{Cognitive Radio Networks: Bernoulli Reward Model}
We first consider a Bernoulli reward model using cognitive radio as
an example of application. There are $M$ secondary users
independently search for idle channels among $N$ channels. In each
slot, the busy/idle state of each channel is drawn from Bernoulli
distribution with unknown mean, \ie
$f(s;\theta)=\theta^s(1-\theta)^{1-s}$ for $s\in\{0,1\}$. When
multiple secondary users choose the same channel for transmission,
they will collide and no one gets reward (collision model 2).

\begin{figure}[h]
\centerline{
\begin{psfrags}
\scalefig{0.7}\epsfbox{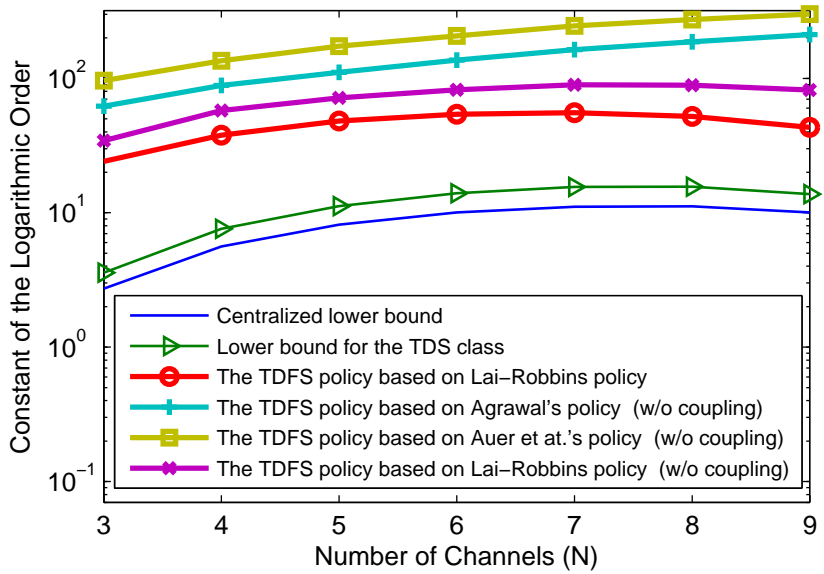}
\end{psfrags}}
\caption{The performance of $\pi^*_F$ built upon different
single-player policies (Bernoulli distributions,~global horizon
length $T=5000,~M=2,~\Theta=[0.1,~0.2,~\cdots,~\frac{N}{10}]$).}
\label{fig:c1}
\end{figure}

In Fig.~\ref{fig:c1}, we plot the performance of $\pi_F^*$ (with
pre-agreement) using different single-player policies. The leading
constant of the logarithmic order is plotted as a function of $N$
with fixed $M$. From Fig.~\ref{fig:c1}, we observe that adopting the
optimal Lai-Robbins single-player policy in $\pi_F^*$ achieves the
best performance. This is expected given that Lai-Robbins policy
achieves the best leading constant in the single-player case. We
also observe that coupling the parallel single-player procedures
through the common observation history in $\pi_F^*$ leads to a
better performance compared to the one without coupling.

\begin{figure}[h]
\centerline{
\begin{psfrags}
\scalefig{0.7}\epsfbox{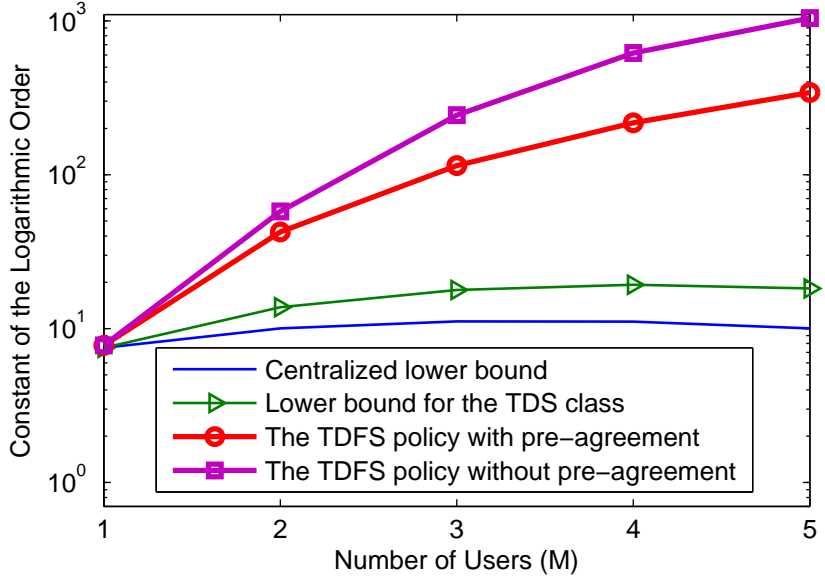}
\end{psfrags}}
\caption{The performance of $\pi^*_F$ based on Lai-Robbins policy
(Bernoulli distributions,~global horizon length
$T=5000,~N=9,~\Theta=[0.1,~0.2,~\cdots,~0.9]$).} \label{fig:u1}
\end{figure}

In Fig.~\ref{fig:u1}, we study the impact of eliminating
pre-agreement on the performance of $\pi_F^*$. We plot the leading
constant of the logarithmic order as a function of $M$ with fixed
$N$.  We observe that eliminating pre-agreement comes with a price
in performance in this case. We also observe that the system
performance degrades as $M$ increases. One potential cause is the
fairness property of the policy that requires each player to learn
the entire rank of the $M$ best arms.

\subsection{Multichannel Communications under Unknown Fading: Exponential Reward Model}

In this example, we consider opportunistic transmission over
wireless fading channels with unknown fading statistics. In each
slot, each user senses the fading condition of a selected channel
and then transmits data with a fixed power. The reward obtained from
a chosen channel is measured by its capacity (maximum data
transmission rate) $C=\log(1+\mbox{SNR})$. We consider the Rayleigh
fading channel model where the SNR of each channel is exponential
distributed with an unknown mean. When multiple users choose the
same channel, no one succeeds (collision model 2). Note that a
channel that has a higher expected SNR also offers a higher expected
channel capacity. It is thus equivalent to consider SNR as the
reward which is exponential distributed. Specifically, we have
$f(s;\theta)=1/\theta \exp(-x/\theta)$ for $s\in\Rc^+$.

\begin{figure}[h]
\centerline{
\begin{psfrags}
\scalefig{0.7}\epsfbox{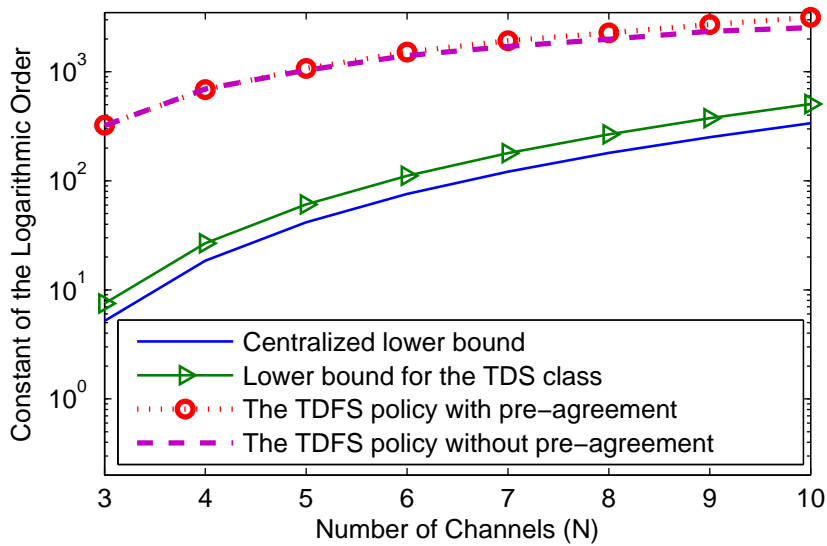}
\end{psfrags}}
\caption{The performance of $\pi^*_F$ based on Agrawal's index
policy (Exponential distributions,~global horizon length
$T=5000,~M=2,~\Theta=[1,~2,~\cdots,~N]$).} \label{fig:c3}
\end{figure}

In Fig.~\ref{fig:c3}, we plot the leading constant of the
logarithmic order as a function of $N$ with fixed $M$. We observe
that in this example, eliminating pre-agreement has little impact on
the performance of the TDFS policy.

\subsection{Target Collecting in Multi-agent Systems: Gaussian Reward Model}
In this example, we consider the Gaussian reward model arisen in the
application of multi-agent systems. At each time, $M$ agents
independently select one of $N$ locations to collect targets (\eg
fishing or ore mining). The reward at each location is determined by
the fish size or ore quality, which has been shown in [5,~6] to fit
with log-Gaussian distribution. The reward at each location has an
unknown mean but a known variance which is the same across all
locations. When multiple agents choose the same location, they share
the reward (collision model 1).

Note that the log-Gaussian and Gaussian reward distributions are
equivalent since the arm ranks under these two distributions are the
same when the arms have the same variance. We can thus focus on the
Gaussian reward model, \ie
$f(s;\theta)=(2\pi\sigma^2)^{-1/2}\exp\{-(s-\theta)^2/(2\sigma^2)\}$
for $s\in\Rc$.

\begin{figure}[h]
\centerline{
\begin{psfrags}
\scalefig{0.7}\epsfbox{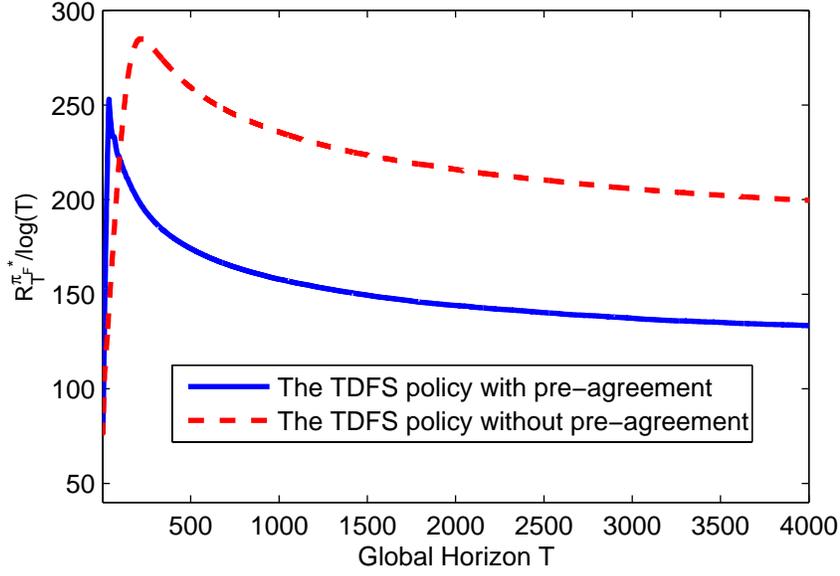}
\end{psfrags}}
\caption{The convergence of the regret growth rate under $\pi^*_F$
based on Lai-Robbins policy (Gaussian
distributions,~$\sigma=1$,~$M=4,~N=9,~\Theta=[1,~2,~\cdots,~9]$).}
\label{fig:r2}
\end{figure}

\begin{figure}[h]
\centerline{
\begin{psfrags}
\scalefig{0.7}\epsfbox{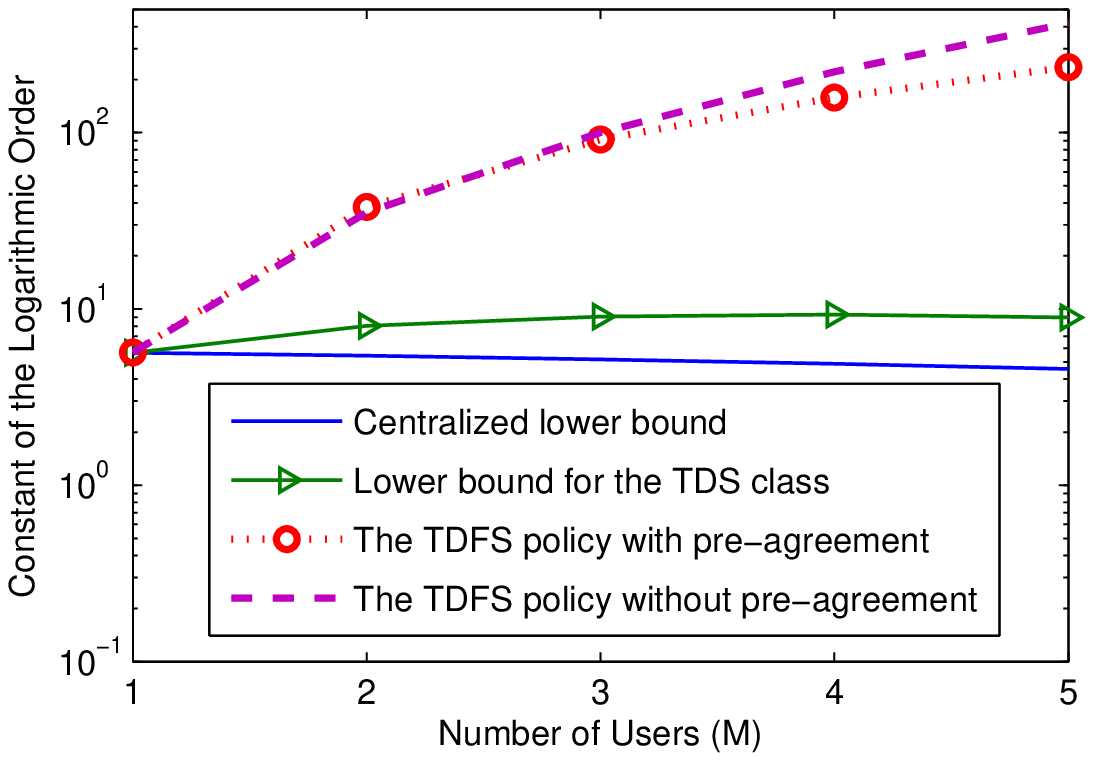}
\end{psfrags}}
\caption{The performance of $\pi^*_F$ based on Lai-Robbins policy
(Gaussian distributions,~$\sigma=1$,~global horizon length
$T=5000,~N=10,~\Theta=[1,~2,~\cdots,~10]$).} \label{fig:u2}
\end{figure}

In Fig.~\ref{fig:r2}, we plot $R_T^{\pi_F^*}/\log(T)$ as a function
of global horizon $T$. We observe that the regret growth rate
quickly converges as time goes, which implies that the policy can
achieve a strong performance within a short finite period.
In Fig.~\ref{fig:u2}, we plot the leading constant of the
logarithmic order as a function of $M$ with fixed $N$. Similar to
the previous example, eliminating pre-agreement has little impact on
the performance of the TDFS policy.

\section{Conclusion} \label{conclusion}
In this paper, we have studied a decentralized formulation of the
classic multi-armed bandit problem by considering multiple
distributed players. We have shown that the optimal system regret in
the decentralized MAB grows at the same logarithmic order as that in
the centralized MAB considered in the classic work by Lai and
Robbins~\cite{Lai&Robbins85AAM} and Anantharam,
\etal\cite{Anantharam&etal87TAC}. A decentralized policy that
achieves the optimal logarithmic order has been constructed. A lower
bound on the leading constant of the logarithmic order is
established for polices in the TDS class, to which the proposed
policy belongs. Future work includes establishing tighter lower
bound on the leading constant of the logarithmic order and
investigating decentralized MAB where each arm has different means
to different players.

\section*{Appendix A.}\label{sec:AppendixA}
Let \textcircled{\footnotesize$H$} denote the set of all possible
values of the unknown parameter $\theta$.

\noindent \textbf{Regularity Conditions:}
\begin{itemize}
\item[C1]Existence of mean: $\mu(\theta)$ exists for any $\theta\in$\textcircled{\footnotesize$H$}.
\item[C2]Positive distance: $0<I(\theta,\theta')<\infty~\mbox{whenever}~~\mu(\theta)<\mu(\theta')$.
\item[C3]Continuity of $I(\theta,\theta')$: $\forall\epsilon>0~\mbox{and}~\mu(\theta)<\mu(\theta'),
\exists\delta>0~\mbox{such that}~
|I(\theta,\theta')-I(\theta,\theta'')|<\epsilon~\mbox{whenever}~
\mu(\theta')\le\mu(\theta'')\le\mu(\theta')+\delta$.
\item[C4]Denseness of \textcircled{\footnotesize$H$}:
$\forall\theta~\mbox{and}~\delta>0,\exists\theta'\in{\it
\Theta}~\mbox{such
that}~\mu(\theta)<\mu(\theta')<\mu(\theta)+\delta$.
\end{itemize}
\textbf{Conditions on point estimate and confidence upper bound:}
{\begin{itemize}
\item[C5] For any $\theta\in$\textcircled{\footnotesize$H$}, we have
{\begin{itemize}\item[(i)]
$\Pr_{\theta}\{g_{t,k}(S(1),\cdots,S(k))\ge r~\mbox{for all}~k<
t\}=1-o(t^{-1})~\mbox{for every}~r<\mu(\theta);$
\item[(ii)]
$
\lim_{\epsilon\downarrow0}(\limsup_{t\rightarrow\infty}\Sigma_{k=1}^{t-1}\Pr_{\theta}\{g_{t,k}(S(1),\cdots,S(k))\ge
\mu(\theta')-\epsilon\}/\log t)\le1/I(\theta,\theta')$\\
$\mbox{whenever}~\mu(\theta')>\mu(\theta);$
\item[(iii)]
$g_{t,k}~\mbox{is nondecreasing in}~t~\mbox{for every
fixed}~k=1,2,\cdots;$ \item[(iv)] $g_{t,k}\ge h_k~\forall~t;$
\item[(v)]
$\Pr_{\theta}\{\max_{\delta t\le k\le
t}|h_k-\mu(\theta)|>\epsilon\}=o(t^{-1})~\forall~\epsilon>0,~0<\delta<1.$
\end{itemize}}
\end{itemize}}

\section*{Appendix B. Proof of
Lemma~\ref{lemma:equiv}}\label{sec:AppendixB} Under the condition
that the $M$ best arms have nonnegative means, the maximum expected
immediate reward obtained under the ideal scenario that $\Theta$ is
known is given by $\Sigma_{j=1}^M\mu(\theta_{\sigma(j)})$. Under a
uniformly good policy $\pi$, we have, for any arm $n$ with
$\mu(\theta_n)\ge\mu(\theta_{\sigma(M)})$,
\[\mathbb{E}[\tau_{n,T}]=T-o(T^b),~\forall~b>0.\]
Following Theorem~3.1 in~\cite{Anantharam&etal87TAC}, the above
equation implies that for any arm $j$ with
$\mu(\theta_j)<\mu(\theta_{\sigma(M)})$,
\[\liminf_{T\rightarrow\infty}\frac{\mathbb{E}[\mathbb{E}[\tau_{j,T}]}{\log T}\ge\frac{1}{I(\theta_j,\theta_{\sigma(M)})}.\]
The regret growth rate is thus lower bounded by
\[\liminf_{T\rightarrow\infty}\frac{R_T(\Theta)}{\log T}\ge\Sigma_{j:\mu(\theta_j)<\mu(\theta_{\sigma(j)})}
\frac{\mu(\theta_{\sigma(M)})-\mu(\theta_j)}{I(\theta_j,\theta_{\sigma(M)})}.\]
Since the optimal policy that chooses exactly $M$ arms at each time
achieves the above lower bound~[1,2], it is also optimal for the MAB
that chooses up to $M$ arms at each time. The two problems are thus
equivalent.

\section*{Appendix C. Proof of Lemma~\ref{lemma:upbadtime}}\label{sec:AppendixC}
Let $l_t$ denote the leader among the arm set
$\{\sigma(i),\cdots,\sigma(N)\}$. Consider arm $n$ with
$\mu(\theta_n)<\mu(\theta_{\sigma(i)})$. Let $\Dc(T)$ denote the set
of slots in the dominant mini-sequence up to time $T$.

For any
$\epsilon\in(0,\mu(\theta_{\sigma(i)})-\mu(\theta_{\sigma(i+1)}))$,
let $N_1(T)$ denote the number of slots in $\Dc(T)$ at which arm $n$
is played when the leader is the $i$th best arm and the difference
between its point estimate and true mean does not exceed $\epsilon$,
$N_2(T)$ the number of slots in $\Dc(T)$ at which arm $n$ is played
when the leader is the $i$th best arm and the difference between its
point estimate and true mean exceeds $\epsilon$, and $N_3(T)$ the
number of slots in $\Dc(T)$ when the leader is not the $i$th best
arm. Recall that each arm is played once in the first $N$ slots. We
have {\begin{eqnarray}\label{eqn:c1} \tau_{n,T}^D\le
1+N_1(T)+N_2(T)+N_3(T).
\end{eqnarray}}Next, we show that $\mathbb{E}[N_1(T)]$, $\mathbb{E}[N_2(T)]$, and $\mathbb{E}[N_3(T)]$
are all at most in the order of $\log T$.

Consider first $\mathbb{E}[N_1(T)]$. Based on the structure of
Lai-Robbins single-player policy, we have {\small\begin{eqnarray}\nn
N_1(T)&=&|\{1<t\le
T:~t\in\{\Dc(T)\},~\mu(\theta_{l_t})=\mu(\theta_{\sigma(i)}),~|h_{\tau_{l_t,t}}-\mu(\theta_{\sigma(i)})|\le\epsilon~\mbox{and
arm}~n~\mbox{is played at time}~t\}|\\\nn&\le& 1+|\{1\le j\le T-1:
g_{s,j}(S_{n1},\cdots,S_{nj})\ge\mu(\theta_{\sigma(i)})-\epsilon~\mbox{for
some}~j<s\le t\}|\\ &\le&1+|\{1\le j\le T-1:
g_{T,j}(S_{n1},\cdots,S_{nj})\ge\mu(\theta_{\sigma(i)})-\epsilon~\}|.
\end{eqnarray}}
Under condition C5, for any $\rho>0$, we can choose $\epsilon$
sufficiently small such that
\begin{eqnarray}\nn
&&\limsup_{T\rightarrow\infty}\frac{\mathbb{E}(|\{1\le j\le T-1:
g_{T,j}(S_{n1},\cdots,S_{nj})\ge\mu(\theta_{\sigma(i)})-\epsilon~\}|)}{\log
T}\\\nn&=& \limsup_{T\rightarrow\infty}\frac{\Sigma_{j=1}^{T-1}\Pr
\{g_{T,j}(S_{n1},\cdots,S_{nj})\ge\mu(\theta_{\sigma(i)})-\epsilon~\}}{\log
T}\\\nn&\le&\frac{1+\rho}{I(\theta_n,\theta_{\sigma(i)})}.
\end{eqnarray}
Thus,
\begin{eqnarray}
\limsup_{T\rightarrow\infty}\frac{\mathbb{E}[N_1(T)]}{\log
T}\le\frac{1}{I(\theta_n,\theta_{\sigma(i)})}.\label{eqn:pc1}
\end{eqnarray}
Consider $\mathbb{E}[N_2(T)]$. Define $c_D(t)\defeq|\{s:s\le
t,s,t\in \Dc(T)\}|$ as the number of slots in $\Dc(T)$ that are no
larger than $t$. Since the number of observations obtained from
$l_t$ is at least $(c_D(t)-1)\delta $, under condition C5, we have,
\begin{eqnarray}\nn
&&\Pr\{c_D(t)=s~\mbox{for some}~s\le t,
\mu(\theta_{l_t})=\mu(\theta_{\sigma(i)}),~|h_{\tau_{l_t,t}}
-\mu(\theta_{\sigma(i)})|>\epsilon\}\\\nn &\le&\Pr\{\sup_{j\ge\delta
(s-1)}|h_j(S_{l_t,1},\cdots,S_{l_t,j})-\mu(\theta_{l_t})|>\epsilon\}\\\nn
&=&\Sigma_{i=0}^{\infty}\delta^io(s^{-1})\\\label{eqn:h}&=&o(s^{-1}).
\end{eqnarray}
We thus have
\begin{eqnarray}\nn
\mathbb{E}[N_2(T)]&=&\mathbb{E}(|\{1<t\le
T:~t\in\Dc(T),~\mu(\theta_{l_t})=\mu(\theta_{\sigma(i)}),
~|h_{\tau_{l_t,t}}-\mu(\theta_{\sigma(i)})|>\epsilon\}|)\\\nn
&\le&\Sigma_{s=1}^{T}\Pr\{c_D(t)=s~\mbox{for some}~s\le t\le T,
\mu(\theta_{l_t})=\mu(\theta_{\sigma(i)}),~|h_{\tau_{l_t,t}}-\mu(\theta_{\sigma(i)})|>\epsilon\}\\\label{eqn:pc2}
&=&o(\log T).
\end{eqnarray}

Next, we show that $\mathbb{E}[N_3]=o(\log T)$.

Choose
$0<\epsilon_1<(\mu(\theta_{\sigma(i)})-\mu(\theta_{\sigma(i+1)}))/2$
and $c>(1-N\delta)^{-1}$. For $r=0,1,\cdots$, define the following
events. {\begin{eqnarray}\nn A_r&\defeq&\cap_{i\le k\le
N}\{\max_{\delta c^{r-1}\le
s}|h_s(S_{\sigma(k),1},\cdots,S_{\sigma(k),s})-\mu(\theta_{\sigma(k)})|\le\epsilon_1\},\\\nn
B_r&\defeq&\{g_{s_m,j}(S_{\sigma(i),1},\cdots,S_{\sigma(i),j})\ge\mu(\theta_{\sigma(i)})-\epsilon~\mbox{for
all}~1\le j\le\delta m,~c^{r-1}\le m\le c^{r+1},~\mbox{and}~s_m>
m\}.
\end{eqnarray}}
By~\eqref{eqn:h}, we have $\Pr(\bar{A_r})=o(c^{-r})$. Consider the
following event:
\begin{eqnarray}
C_r\defeq\{g_{m+1,j}(S_{\sigma(i),1},\cdots,S_{\sigma(i),j})\ge\mu(\theta_{\sigma(i)})-\epsilon_1~\mbox{for
all}~1\le i\le\delta m,~\mbox{and}~c^{r-1}\le m\le c^{r+1}\}.
\end{eqnarray}
Under condition C5, it is easy to see that $B_r\supset C_r$. From
Lemma $1-(i)$ in~\cite{Lai&Robbins85AAM},
$\Pr(\bar{C_r})=o(c^{-r})$. We thus have $\Pr(\bar{B_r})=o(c^{-r})$.

Consider time $t$. Define the event $D_r\defeq\{t\in
\Dc(T),~c^{r-1}\le c_D(t)-1<c^{r+1}\}$. When the round-robin
candidate $r_t=\sigma(i)$, we show that on the event $A_r\cap
B_r\cap D_r$, $\sigma(i)$ must be played. It is sufficient to focus
on the nontrivial case that
$\mu(\theta_{l_t})<\mu(\theta_{\sigma(i)})$. Since $\tau_{l_t,t}\ge
(c_D(t)-1)\delta$, on $A_r\cap D_r$, we have
$h_{\tau_{l_t,t}}<\mu(\theta_{\sigma(i)})-\epsilon_1$. We also have,
on $A_r\cap B_r\cap D_r$,
\begin{eqnarray}
g_{t,\tau_{\sigma(i),t}}(S_{\sigma(i),1},\cdots,S_{\sigma(i),\tau_{\sigma(i),t}})\ge\mu(\theta_{\sigma(i)})-\epsilon_1.
\end{eqnarray}
Arm $\sigma(i)$ is thus played on $A_r\cap B_r\cap D_r$. Since
$(1-c^{-1})/N>\delta$, for any $c^{r}\le s\le c^{r+1}$, there exists
an $r_0$ such that on $A_r\cap B_r$, $\tau_{\sigma(i),t}\ge
(1/N)(s-c^{r-1}-2N)>\delta s$ for all $r>r_0$. It thus follows that
on $A_r\cap B_r\cap D_r$, for any $c^{r}\le c_D(t)-1\le c^{r+1}$, we
have $\tau_{\sigma(i),t}>(c_D(t)-1)\delta$, and $\sigma(i)$ is thus
the leader. We have, for all $r>r_0$,
\begin{eqnarray}
\Pr(l_t\neq\sigma(i),t\in \Dc(T),c^r\le c_D(t)-1< c^{r+1})\le
\Pr(\bar{A_r}\cap D_r)+\Pr(\bar{B_r}\cap D_r)=o(c^{-r}).
\end{eqnarray}
Therefore, {\begin{eqnarray}\nn
\mathbb{E}[N_3(T)]&=&\mathbb{E}[|\{t>1:~t\in\Dc(T),~h_{\tau_{l_t,t}}<\mu(\theta_{\sigma(i)})\}|]\\\nn
&\le&\Sigma_{s=0}^{T-1}\Pr(c_D(t)-1=s~\mbox{for
some}~t>s,l_t\neq\sigma(i), t\in \Dc(T))\\\nn&\le&
1+\Sigma_{r=0}^{\lceil\log_cT\rceil}\Sigma_{c^r\le s\le
c^{r+1}}\Pr(\exists t,~l_t\neq\sigma(i),t\in \Dc(T),
c_D(t)-1=s)\\\nn
&=&\Sigma_{r=0}^{\lceil\log_cT\rceil}o(1)\\\label{eqn:pc3}&=& o(\log
T).
\end{eqnarray}}
From~\eqref{eqn:pc1},~\eqref{eqn:pc2},~\eqref{eqn:pc3}, we arrive at
Lemma~\ref{lemma:upbadtime}.

\section*{Appendix D. Proof of
Lemma~\ref{lemma:lbbesttime}}\label{sec:AppendixD} We prove by
induction on $i~(1\le i\le M)$. Consider $i=1$. In this case,
Lai-Robbins policy is applied to all arms by player 1. Let
$\mathbb{E}[\hat{\tau}_{\overline{\sigma(i)},T}]$ denote the
expected total number of times in the first sequence that player $i$
does not play the $i$th best arm up to time $T$. From
Lemma~\ref{lemma:upbadtime}, we have
\begin{eqnarray}
\limsup_{T\rightarrow\infty}\frac{\mathbb{E}[\hat{\tau}_{\overline{\sigma(1)},T}]}{\log
T}=\Sigma_{n:\mu(\theta_n)<\mu(\theta_{\sigma(1)})}\limsup_{T\rightarrow\infty}\frac{\mathbb{E}[\tau^D_{n,T}]}{\log
T}\le x_1.
\end{eqnarray}
Since
$\frac{(T-2M)}{M}\le\mathbb{E}[\hat{\tau}_{\overline{\sigma(1)},T}+\hat{\tau}_{\sigma(1),T}]\le
\frac{T}{M}$, we have
\begin{eqnarray}
\limsup_{T\rightarrow\infty}\frac{T/M-\mathbb{E}[\hat{\tau}_{\sigma(1),T}]}{\log
T}=\limsup_{T\rightarrow\infty}\frac{\mathbb{E}[\hat{\tau}_{\overline{\sigma(1)},T}]}{\log
T}\le x_1,
\end{eqnarray}
which shows that~\eqref{eqn:lb1} holds for $i=1$. Let
$\mathbb{E}[\tau^{D_i}_{T}]$ denote the number of slots in the
dominant mini-sequence of player $i$ up to time $T$.
From~\eqref{eqn:lb1}, we have, for $i=1$,
\begin{eqnarray}\label{eqn:lb2}
\limsup_{T\rightarrow\infty}\frac{T/M-\mathbb{E}[\tau^{D_{i+1}}_{T}]}{\log
T}\le x_{i}.
\end{eqnarray}

Assume~\eqref{eqn:lb1} and~\eqref{eqn:lb2} hold for $i=k$. To
complete the induction process, consider $i=k+1$. Let
$\mathbb{E}[\tau^{D_i}_{\sigma(i),T}]$ denote the expected number of
slots that the $i$th best arm is played by player $i$ in its
dominant mini-sequence up to time $T$, and
$\mathbb{E}[\tau^{D_i}_{\overline{\sigma(i)},T}]$ the expected
number of slots that the $i$th best arm is not played by player $i$
in its dominant mini-sequence up to time $T$. Based on
Lemma~\ref{lemma:upbadtime}, we have
\begin{eqnarray}
\limsup_{T\rightarrow\infty}\frac{\mathbb{E}[\tau^{D_i}_{\overline{\sigma(i)},T}]}{\log
T}\le\Sigma_{n:\mu(\theta_n)<\mu(\theta_{\sigma(i)})}\frac{1}{I(\theta_n,\theta_{\sigma(i)})}.
\end{eqnarray} Since
$\mathbb{E}[\tau^{D_i}_{T}]=\mathbb{E}[\tau^{D_i}_{\overline{\sigma(i)},T}+\tau^{D_i}_{\sigma(i),T}]
\le\mathbb{E}[\tau^{D_i}_{\overline{\sigma(i)},T}+\hat{\tau}_{\sigma(i),T}]$,
we have
\begin{eqnarray}\nn
\limsup_{T\rightarrow\infty}\frac{T/M-\mathbb{E}[\hat{\tau}_{\sigma(i),T}]}{\log
T}&\le&\limsup_{T\rightarrow\infty}\frac{T/M-\mathbb{E}[\tau^{D_i}_{\sigma(i),T}]}{\log
T}\\\nn&=&\limsup_{T\rightarrow\infty}\frac{T/M-\mathbb{E}[\tau^{D_i}_{T}]+\mathbb{E}[\tau^{D_i}_{\overline{\sigma(i)},T}]}{\log
T}\\&\le& x_{k+1}.
\end{eqnarray}
We thus proved Lemma~\ref{lemma:lbbesttime} by induction.

\section*{Appendix E. Proof of
the Upper Bound $C(\Theta)$ in
Theorem~\ref{thm:worstcase}}\label{sec:AppendixE}

\noindent{\em Collision Model 1.}Consider the first subsequence.
Since the expected number of slots that the $i$th best arm is played
by player $i$ dominates that by other players, the expected total
reward $\mathbb{E}[Y_{i}(T)]$ obtained from the $i$th best arm by
all players up to time $T$ satisfies
\begin{eqnarray}\label{eqn:r2}
\limsup_{T\rightarrow\infty}\frac{T\mu(\theta_{\sigma(i)})/M-\mathbb{E}[Y_{i}(T)]}{\log
T}\le x_{i}\mu(\theta_{\sigma(i)}).
\end{eqnarray}
For an arm $n$ with $\mu(\theta_n)<\mu(\theta_{\sigma(M)})$, the
worst case is that player $M$ plays arm $n$ for
$1/I(\theta_n,\theta_{\sigma(M)})\log T-o(\log T)$ times up to time
$T$ (by Lemma~\ref{lemma:upbadtime}). From~\eqref{eqn:r2}, the
expected total reward $\Sigma_{i=1}^N\mathbb{E}[Y_{i}(T)]$ obtained
from all arms and players up to time $T$ satisfies
\begin{eqnarray}\nn
&&\limsup_{T\rightarrow\infty}\frac{T\Sigma_{i=1}^M\mu(\theta_{\sigma(i)})/M-\Sigma_{i=1}^N\mathbb{E}[Y_{i}(T)]}{\log
T}\\&\le&
\Sigma_{i=1}^Mx_{i}\mu(\theta_{\sigma(i)})-\Sigma_{n:\mu(\theta_n)<\mu(\theta_{\sigma(M)})}\frac{1}{I(\theta_n,\theta_{\sigma(M)})}\mu(\theta_n).\label{eqn:r3}
\end{eqnarray}
The expected total reward
$\Sigma_{j=1}^M\Sigma_{i=1}^N\mathbb{E}[Y_{i}(T)]$ obtained from all
arms and players over all the $M$ subsequences up to time $T$ thus
satisfies
\begin{eqnarray}\nn
&&\limsup_{T\rightarrow\infty}\frac{T\Sigma_{i=1}^M\mu(\theta_{\sigma(i)})-\Sigma_{j=1}^M\Sigma_{i=1}^N\mathbb{E}[Y_{i}(T)]}{\log
T}\\&\le&
M(\Sigma_{i=1}^Mx_{i}\mu(\theta_{\sigma(i)})-\Sigma_{n:\mu(\theta_n)<\mu(\theta_{\sigma(M)})}\frac{1}{I(\theta_n,\theta_{\sigma(M)})}\mu(\theta_n)).\label{eqn:r4}
\end{eqnarray}

\noindent{\em Collision Model 2.} Consider the first subsequence. By
Lemma~\ref{lemma:lbbesttime}, the expected number of slots in which
the $i$th best arm is played by player $j~(j\neq i)$ that is at most
$x_j\log T+o(\log T)$ up to time $T$. Together with~\eqref{eqn:lb1},
the expected total reward $\mathbb{E}[Y_{i}(T)]$ obtained from arm
$i$ by all players up to time $T$ satisfies
\begin{eqnarray}\label{eqn:r5}
\limsup_{T\rightarrow\infty}\frac{T\mu(\theta_{\sigma(i)})/M-\mathbb{E}[Y_{i}(T)]}{\log
T}\le \Sigma_{k=1}^Mx_{k}\mu(\theta_{\sigma(i)}).
\end{eqnarray}
For an arm $n$ with $\mu(\theta_n)<\mu(\theta_{\sigma(M)})$, the
worst case is that player $i$ plays arm $n$ for
$1/I(\theta_n,\theta_{\sigma(i)})\log T-o(\log T)$ times up to time
$T$. Combined with the worst case collisions, the expected reward
obtained on arm $n$ is at least
$(\max\{1/I(\theta_n,\theta_M)-\Sigma_{i=1}^{M-1}1/I(\theta_n,\theta_{\sigma(i)}),0\}\mu(\theta_n))\log
T-o(\log T)$ up to time $T$. From~\eqref{eqn:r5}, the expected total
reward $\Sigma_{i=1}^N\mathbb{E}[Y_{i}(T)]$ obtained from all arms
and players up to time $T$ satisfies {\begin{eqnarray}\nn
&&\limsup_{T\rightarrow\infty}\frac{T\Sigma_{i=1}^M\mu(\theta_{\sigma(i)})/M-\Sigma_{i=1}^N\mathbb{E}[Y_{i}(T)]}{\log
T}\\&\le&
\Sigma_{i=1}^M\Sigma_{k=1}^Mx_{k}\mu(\theta_{\sigma(i)})-\Sigma_{n:\mu(\theta_n)<\mu(\theta_{\sigma(M)})}
\max\{\frac{1}{I(\theta_n,\theta_M)}-\Sigma_{i=1}^{M-1}\frac{1}{I(\theta_n,\theta_{\sigma(i)})},0\}\mu(\theta_n).~~~~~\label{eqn:r6}
\end{eqnarray}}
The expected total reward
$\Sigma_{j=1}^M\Sigma_{i=1}^N\mathbb{E}[Y_{i}(T)]$ obtained from all
arms and players over all the $M$ subsequences up to time $T$ thus
satisfies
\begin{eqnarray}\nn
&&\limsup_{T\rightarrow\infty}\frac{T\Sigma_{i=1}^M\mu(\theta_{\sigma(i)})-\Sigma_{j=1}^M\Sigma_{i=1}^N\mathbb{E}[Y_{i}(T)]}{\log
T}\\&\le&
M(\Sigma_{i=1}^M\Sigma_{k=1}^Mx_{k}\mu(\theta_{\sigma(i)})-\Sigma_{n:\mu(\theta_n)<\mu(\theta_{\sigma(M)})}
\max\{\frac{1}{I(\theta_n,\theta_M)}-\Sigma_{i=1}^{M-1}\frac{1}{I(\theta_n,\theta_{\sigma(i)})},0\}\mu(\theta_n)).~~~~~\label{eqn:r7}
\end{eqnarray}

\section*{Appendix F. Proof of
Theorem~\ref{thm:optElimPre}}\label{sec:AppendixF}

Based on the structure of $\pi_F^*$ without pre-agreement, each
player's local policy consists of disjoint rounds, each of which
consists of $M$ slots where $M$ arms are played with a particular
local offset. We say that a round is in group $o~(0\le o\le M-1)$ if
it corresponds to local offset $o$. We also say that a round is {\em
normal} if the player plays the $M$ best arms in the correct order
according to its local offset, otherwise, the round is {\em
singular}. Let $t$ denote the global reference time. We say that
slot $t$ is {\em normal} if every player's round that contains slot
$t$ is normal, otherwise, slot $t$ is {\em singular}.

Note that the set of slots where a reward loss incurs is a subset of
the singular slots and normal slots involving collisions. To prove
the logarithmic order of the system regret growth rate, it is thus
sufficient to show the following: (i) the expected number of
singular slots is at most logarithmic with time; (ii) the expected
number of normal slots involving collisions is at most logarithmic
with time.

To show (i), we consider an arbitrary player. Similar to the proof
of Theorem~\ref{thm:worstcase}, it can be shown that in any group
$o$, the expected number of singular rounds is at most logarithmic
with time. Summed over all $M$ groups, the total expected number of
singular rounds is at most logarithmic with time. Since that a
singular round corresponds to $M$ singular slots, the expected
number of singular slots is at most logarithmic with time.

To show (ii), we first show that between any two consecutive
singular slots, the expected number of normal slots involving
collisions is uniformly bounded by some constant (independent of the
two singular slots being considered). Note that a collision occurs
in a normal slot only when multiple players have the same global
offset (with respect to the global time $t$). Since a player will
randomly generate a new local offset if a collision is observed, it
is easy to show that, between any two consecutive singular slots,
the expected number of normal slots needed for all players to settle
at different global offsets is uniformly bounded. Consequently,
between any two consecutive singular slots, the expected number of
normal slots involving collisions is uniformly bounded. It thus
follows that the expected number of normal slots involving
collisions is at the same logarithmic order as the expected number
of singular slots.

The fairness of $\pi_F^*$ without pre-agreement follows directly
from the symmetry among players.

\section*{Appendix G. Proof of
Lemma~\ref{lemma:lowboundgen}}\label{sec:AppendixG} For any
$\delta>0$, under condition C1-C4, we can choose a parameter
$\lambda$ such that
\[\mu(\theta_{\sigma(j-1)})>\mu(\lambda)>\mu(\theta_{\sigma(j)})>\mu(\theta_{\sigma(k)})\]
and
\begin{eqnarray}\label{eqn:l0}
|I(\theta_{\sigma(k)},\lambda)-I(\theta_{\sigma(k)},\theta_{\sigma(j)})|\le\delta
I(\theta_{\sigma(k)},\theta_{\sigma(j)}).
\end{eqnarray}Consider the new
parameter set $\Theta'$ where the mean of the $k$th best arm (say
arm $l$) is replaced by $\lambda$. Since arm $l$ is the $j$th best
arm under $\Theta'$, from~\eqref{eqn:lowtime}, we have, for any
$0<b<\delta$,
\begin{eqnarray}
\mathbb{E}_{\Theta'}[v(T)-\tau_{l,T}]\le o((v(T))^b),
\end{eqnarray}
where $\tau_{l,T}$ is the number of times that arm $l$ is played up
to time $T$. There thus exists a random function $c(T)$ with
$c(T)+v(T)-\tau_{l,T}\ge 0$ a.s. and $\mathbb{E}_{\Theta'}[c(T)]\ge
0$ such that\footnote{For example, $c(T)$ can be constructed as
$c(T)=b(T)-v(T)+\tau_{l,T}$, where
$b(T)=[\mathbb{E}_{\Theta'}[v(T)-\tau_{l,T}]]^+.$}
\begin{eqnarray}\label{eqn:ll}
\mathbb{E}_{\Theta'}[c(T)+v(T)-\tau_{l,T}]=o((v(T))^b).
\end{eqnarray}
We have {\small\begin{eqnarray}\nn
\mathbb{E}_{\Theta'}[c(T)+v(T)-\tau_{l,T}]&=&\mathbb{E}_{\Theta'}[c(T)+v(T)-\tau_{l,T}|\tau_{l,T}<\frac{(1-\delta)\log
(v(T))}{I(\theta_{\sigma(k)},\lambda)}]\textrm{Pr}_{\Theta'}\{\tau_{l,T}<\frac{(1-\delta)\log
(v(T))}{I(\theta_{\sigma(k)},\lambda)}\}\\\nn&&+\mathbb{E}_{\Theta'}[c(T)+v(T)-\tau_{l,T}|\tau_{l,T}\ge\frac{(1-\delta)\log
(v(T))}{I(\theta_{\sigma(k)},\lambda)}]\textrm{Pr}_{\Theta'}\{\tau_{l,T}\ge\frac{(1-\delta)\log
(v(T))}{I(\theta_{\sigma(k)},\lambda)}\}\\\nn&\ge&\mathbb{E}_{\Theta'}[c(T)+v(T)-\tau_{l,T}|\tau_{l,T}<\frac{(1-\delta)\log
(v(T))}{I(\theta_{\sigma(k)},\lambda)}]\textrm{Pr}_{\Theta'}\{\tau_{l,T}<\frac{(1-\delta)\log
(v(T))}{I(\theta_{\sigma(k)},\lambda)}\}\\\label{eqn:l1}&\ge&(\mathbb{E}_{\theta'}[c(T)]+v(T)-\frac{(1-\delta)\log
(v(T))}{I(\theta_{\sigma(k)},\lambda)})\textrm{Pr}_{\Theta'}\{\tau_{l,T}<\frac{(1-\delta)\log
(v(T))}{I(\theta_{\sigma(k)},\lambda)}\},
\end{eqnarray}}where the fist inequality is due to the fact that
$c(T)+v(T)-\tau_{l,T}\ge0$ almost surely.

Let $S_{l,1},S_{l,2},\cdots$ denote independent observations from
arm $l$. Define
\[L_t\defeq\Sigma_{i=1}^t\log(\frac{f(S_{l,i};\theta_{\sigma(k)})}{f(S_{l,i};\lambda)})\]
and event \[C\defeq\{\tau_{l,T}<\frac{(1-\delta)\log
(v(T))}{I(\theta_{\sigma(k)},\lambda)},~L_{\tau_{l,T}}\le (1-b)\log
v(T)\}.\] Following~\eqref{eqn:ll} and~\eqref{eqn:l1}, we have
{\small\begin{eqnarray}\nn
\textrm{Pr}_{\Theta'}\{C\}&\le&\textrm{Pr}_{\Theta'}\{\tau_{l,T}<\frac{(1-\delta)\log
(v(T))}{I(\theta_{\sigma(k)},\lambda)}\}\\\nn &\le&
\frac{o((v(T))^b)}{(\mathbb{E}_{\theta'}[c(T)]+v(T)-(1-\delta)\log
(v(T))/I(\theta_{\sigma(k)},\lambda))}\\\nn&\le&
\frac{o((v(T))^b)}{(v(T)-(1-\delta)\log
(v(T))I(\theta_{\sigma(k)},\lambda))}\\&=&o((v(T))^{b-1}),\label{eqn:l2}
\end{eqnarray}}where the last inequality is due to the fact that $\mathbb{E}_{\Theta'}[c(T)]\ge
0$.

We write $C$ as the union of mutually exclusive events
$C_s\defeq\{\tau_{l,T}=s,~L_{s}\le (1-b)\log v(T)\}$ for each
integer $s<(1-\delta)\log (v(T))/I(\theta_{\sigma(k)},\lambda)$.
Note that
\begin{eqnarray}\nn
\textrm{Pr}_{\theta'}\{C_s\}&=&\int_{\{\tau_{l,T}=s,~L_{s}\le
(1-b)\log
v(T)\}}d\textrm{Pr}_{\theta'}\\\nn&=&\int_{\{\tau_{l,T}=s,~L_{s}\le
(1-b)\log v(T)\}}\prod_{i=1}^s
\frac{f(S_{l,i};\lambda)}{f(S_{l,i};\theta_{\sigma(k)})}d\textrm{Pr}_{\theta}\\
&\ge&\exp(-(1-b)\log v(T))\textrm{Pr}_{\theta}\{C_s\}.\label{eqn:l3}
\end{eqnarray}

Based on~\eqref{eqn:l2} and~\eqref{eqn:l3}, we have
\begin{eqnarray}\label{eqn:l4}
\textrm{Pr}_{\theta}\{C\}=(\log
(v(T)))^{1-b}\textrm{Pr}_{\theta'}\{C\}\rightarrow 0.
\end{eqnarray}

By the strong law of large numbers, $L_t/t\rightarrow
I(\theta_{\sigma(k)},\lambda)>0$ as $t\rightarrow\infty$ a.s. under
$\Pr_{\Theta}$. This leads to $\max_{i\le t}L_i/t\rightarrow
I(\theta_{\sigma(k)},\lambda)>0$ a.s. under $\Pr_{\Theta}$. Since
$1-b>1-\delta$, we have
\begin{eqnarray}
\lim_{T\rightarrow\infty}\textrm{Pr}_{\theta}\{L_i>(1-b)\log
v(T)~\mbox{for some}~i<\frac{(1-\delta)\log
v(T)}{I(\theta_{\sigma(k)},\lambda)}\}=0,
\end{eqnarray}which leads to
\begin{eqnarray}\label{eqn:l5}
\lim_{T\rightarrow\infty}\textrm{Pr}_{\theta}\{\tau_{l,T}<\frac{(1-\delta)\log
(v(T))}{I(\theta_{\sigma(k)},\lambda)},~L_{\tau_{l,T}}>(1-b)\log
v(T)\}=0.
\end{eqnarray}
Based on~\eqref{eqn:l4} and~\eqref{eqn:l5}, we have
\begin{eqnarray}
\lim_{T\rightarrow\infty}\textrm{Pr}_{\theta}\{\tau_{l,T}<\frac{(1-\delta)\log
v(T)}{I(\theta_{\sigma(k)},\lambda)}\}=0.
\end{eqnarray}
Based on~\eqref{eqn:l0}, we arrive at
\begin{eqnarray}
\lim_{T\rightarrow\infty}\textrm{Pr}_{\theta}\{\tau_{l,T}<\frac{(1-\delta)\log
v(T)}{((1+\delta)I(\theta_{\sigma(k)},\theta_{\sigma(j)}))}\}=0~\mbox{for
any}~\delta>0,
\end{eqnarray}which leads to~\eqref{eqn:l}.

\bibliographystyle{ieeetr}
{

}

\end{document}

%% file: macros.tex
\def\boxit#1{\vbox{\hrule\hbox{\vrule\kern3pt
        \vbox{\kern3pt#1\kern3pt}\kern3pt\vrule}\hrule}}

\def\reals{ { {\rm  I \kern-0.15em R }  } }
\def\complex{ {\,{{\rm C} \kern-0.50em \raise0.20ex {  |}}\, }}

\def\Rbf{{\bf R}}

\def\Ac{{\cal A}}
\def\Bc{{\cal B}}

\def\Dc{{\cal D}}

\def\Nc{{\cal N}}

\def\Rc{{\cal R}}

\def\be{\begin{equation}}
\def\ee{\end{equation}}

\def\defeq{{\stackrel{\Delta}{=}}}
\def\scalefig#1{\epsfxsize #1\textwidth}
\def\Rxx{\Rbf_{\ssstyle X\kern-.1em X}}

\let\ssstyle=\scriptscriptstyle

\def\eg{{\it e.g.,\ \/}}
\def\etal{{\it et al. \/}}

\def\ie{{\it i.e.,\ \/}}
\def\Kout{\setbox1=\hbox{\Huge\bf K}\hbox to
1.05\wd1{\hspace{.05\wd1}
}}
\def\Sout{\setbox1=\hbox{\Huge\bf S}\hbox to 1.05\wd1{\hspace{.05\wd1}
}}